\documentclass[11pt]{article}

\usepackage{amsmath}
\usepackage{amssymb}
\usepackage{amsthm}
\usepackage{mathtools}
\usepackage{tabu}
\usepackage{textcomp}
\usepackage{listings}
\usepackage{graphicx}
\usepackage{tikz}
\usepackage{arydshln}
\usepackage[margin=1in]{geometry}
\usepackage{indentfirst}

\usepackage{array}
\usepackage{xltabular}
\usepackage{caption}

\newtheorem{proposition}{Proposition}[section]
\newtheorem{example}{Example}[section]
\newtheorem{problem}{Problem}[section]

\newcommand{\bbN}{\mathbb{N}}
\newcommand{\bbZ}{\mathbb{Z}}

\newcommand{\bbK}{\mathbb{K}}

\newcommand{\rtri}{_{\triangleright}}
\newcommand{\ltri}{{}_{\triangleleft}}

\newcommand{\tikzpic}[2]{
    \begin{tikzpicture}[scale=#1]
        #2
    \end{tikzpicture}}

\newcommand{\drawtile}[2]{
    \draw (#1,#2) rectangle (#1+1,#2+1)}

\newcommand{\drawdomino}[4]{
    \drawtile{#1}{#2};
    \drawtile{#1+#3}{#2+#4}}

\newcommand{\drawtromino}[6]{
    \drawtile{#1}{#2};
    \drawtile{#1+#3}{#2+#4};
    \drawtile{#1+#5}{#2+#6}}

\newcommand{\drawtetromino}[8]{
    \drawtile{#1}{#2};
    \drawtile{#1+#3}{#2+#4};
    \drawtile{#1+#5}{#2+#6};
    \drawtile{#1+#7}{#2+#8}}

\newcommand{\drawwtile}[6]{
    \draw (#1,#2) rectangle (#1+1,#2+1);
    \node at(#1+0.25,#2+0.5) [left] {#3};
    \node at(#1+0.5,#2+0.75) [above] {\phantom{$d$} #4 \phantom{$p$}};
    \node at(#1+0.5,#2+0.25) [below] {\phantom{$d$} #5 \phantom{$p$}};
    \node at(#1+0.75,#2+0.5) [right] {#6}
}

\newcommand{\wtile}[4]{\hspace{-1em}
    \raisebox{-0.3cm-1em}{
        \tikzpic{0.6}{
            \drawwtile{0}{0}{#1}{#2}{#3}{#4};
        }
    }\hspace{-0.6em}
}

\newcommand{\whdomino}[8]{\hspace{-1em}
    \raisebox{-0.3cm-1em}{
        \tikzpic{0.6}{
            \drawwtile{0}{0}{#1}{#2}{#3}{#4};
            \drawwtile{1}{0}{#5}{#6}{#7}{#8};
        }
    }\hspace{-0.6em}
}

\newcommand{\Ii}{\tikzpic{0.15}{\drawtile{0}{0};}}
\newcommand{\Iii}{\tikzpic{0.15}{\drawtile{0}{0};\drawtile{0}{1};}}
\newcommand{\Iiii}{\tikzpic{0.15}{\drawtile{0}{0};\drawtile{0}{1};\drawtile{0}{2};}}
\newcommand{\Liii}{\tikzpic{0.15}{\drawtile{0}{0};\drawtile{0}{1};\drawtile{1}{0};}}
\newcommand{\Oiv}{\tikzpic{0.15}{\drawtile{0}{0};\drawtile{0}{1};\drawtile{1}{0};\drawtile{1}{1};}}

\title{Solving tiling enumeration problems by tensor network contractions}
\author{Kai Liang}
\date{\today} 

\begin{document}
    \maketitle

    \begin{abstract}
        This paper presents an algorithm for computing the contraction of two-dimensional tensor networks on a square lattice;
        and we combine it with solving congruence equations to compute the exact enumeration (including weighted enumeration) of Wang tilings.
        Based on this, the paper demonstrates how to transform other tiling enumeration problems (such as those of polyominoes) into Wang tiling enumeration problems, thereby solving them using this algorithm.

        Our algorithm extends the sequence length records for dozens of sequences defined by polyomino tiling enumeration on chessboards on the OEIS website,
        covering numerous of different polyomino sets, including I-polyominoes, tetrominoes, pentominoes, etc.
        This demonstrates the high efficiency and strong universality of the algorithm for solving exact tiling enumeration problems.

        In addition, the theory and techniques used in the algorithm establish a bridge between tensor network contractions and tiling enumeration, where the former provides a theoretical foundation for solving problems in the latter, while the latter offers an intuitive combinatorial interpretation of the former.
    \end{abstract}
    
    \section{Introduction}
    
    Tiling enumeration is a classic problem in combinatorics~\cite{nOT_prob, til_enum, til_enum_1}.

    Taking the most classic \textit{domino tiling} (\cite{DT, DT_0, DT_1, DT_2, DT_3}) as an example: 
    given a \textit{rectangle chessboard} of size $l \times h$ (i.e., with width $l$ and height $h$),
    the task is to calculate the number of all possible ways to cover the entire chessboard by tiling \textit{dominoes}
    (as shown below, $1 \times 2$ tiles that can be rotated), distinguishing between distinct ways that are symmetric.

    \begin{center}
        \tikzpic{0.6}{
            \drawdomino{0}{0}{0}{1};
            \drawdomino{2}{0}{1}{0};
        }
    \end{center}

    For example, the following is one of the tilings using dominoes in a $4 \times 3$ chessboard:

    \begin{center}
        \tikzpic{0.6}{
            \draw[color=gray, dashed] (0,0) grid (4,3);
            \draw (0,0) rectangle (4,3);
            \draw (0,2)--(4,2);\draw (2,1)--(4,1);\draw (1,0)--(1,2);\draw (2,0)--(2,3);
        }
    \end{center}

    When the shape of the chessboard is $l \times h$, the number of different ways of tiling (also reffered to as \textit{tiling enumeration}) can be calculated using the following formula~\cite{DT, DT_0}:
    \[
        \left(\prod_{i=1}^l \prod_{j=1}^h \left( 4\cos^2\frac{i\pi}{l+1}+4\cos^2\frac{j\pi}{h+1} \right)\right) ^{\frac{1}{4}}.
    \]
    However, most other tiling enumeration problems do not have specific formulas applicable to the general case like this,
    with at most formulas for special cases where the width $l$ is fixed and small~\cite{nOT_form, DT_4}.
    In most cases, such enumeration problems requires complex computation,
    and the time and space complexity of the algorithms are both exponential~\cite{WT_CC, WT_CC_1}.

    The algorithm presented in this paper is designed for the exact enumeration of Wang tilings on various shapes of chessboards, and can be generalized to other tiling problems (such as polyomino tilings).
    Compared to previous algorithms, it effectively saves space and improves efficiency.
    The website OEIS (The On-Line Encyclopedia of Integer Sequences) records hundreds of sequences defined by similar tiling enumeration in chessboards,
    Our algorithm has extended the known sequence lengths for dozens of these cases beyond previous record.

    Section 2 of this paper introduces the preliminary knowledge involved in the algorithm,
    including tiling enumeration and tensor network contractions.

    Section 3 explains the specific principles of the algorithm for the exact enumeration of Wang tilings and analyzes its related parameters.

    Section 4 discusses how to transform domino tilings into Wang tilings, thereby enabling exact enumeration using this algorithm.

    Section 5 provides a summary of some further generalizations of the algorithm, and enumerates several unresolved open problems.

    \section{Preliminary}

    \subsection{Tiling enumeration}

    \textit{Wang tiles} (\cite{WT, WT_1}) are unit square tiles, each with a character on its four edges, such as
    \[
        \wtile{$\alpha$}{$a$}{$b$}{$\beta$}.
    \]
    In general, the characters on the left and right edges of a tile are taken from a character set, which we call the \textit{horizontal character set};
    while the characters on the top and bottom edges are taken from another character set, which we call the \textit{vertical character set}.
    Wang tiles can be joined edge-to-edge, provided that the characters on adjacent edges are identical.
    For example, the horizontal concatenation of two Wang tiles is
    \[
        \wtile{$\alpha_1$}{$a_1$}{$b_1$}{$\beta_1$}\circ\wtile{$\alpha_2$}{$a_2$}{$b_2$}{$\beta_2$}
        =\begin{cases}
             \whdomino{$\alpha_1$}{$a_1$}{$b_1$}{$\beta_1$}{}{$a_2$}{$b_2$}{$\beta_2$}, &\beta_1=\alpha_2;\\
             \emptyset, & \beta_1\neq \alpha_2.
        \end{cases}
    \]

    A \textit{Wang tiling} requires covering a specific region with Wang tiles that match edge-to-edge.
    Under the standard rule, Wang tiles cannot be rotated or flipped,
    and each position on the edge of the given region is assigned a fixed character (called the \textit{boundary character}, usually denoted as $\sharp$),
    requiring the tiles on the edge to match these characters.
    Various variants of Wang tilings in finite regions can be reduced to this standard form.

    For instance, we select the Wang tile set as
    \[
        \{\wtile{$\sharp$}{$\sharp$}{$1$}{$\sharp$},\wtile{$\sharp$}{$\sharp$}{$\sharp$}{$1$},\wtile{$\sharp$}{$1$}{$\sharp$}{$\sharp$},\wtile{$1$}{$\sharp$}{$\sharp$}{$\sharp$}\}.
    \]
    Then, one of the tilings in a $4\times 3$ chessboard is
    \begin{center}
        \raisebox{-1.2cm}{\tikzpic{0.6}{
            \draw (0,0) rectangle (4,3);
            \node at(0.25,0.5) [left] {$\sharp$};\node at(0.25,1.5) [left] {$\sharp$};\node at(0.25,2.5) [left] {$\sharp$};
            \node at(3.75,0.5) [right] {$\sharp$};\node at(3.75,1.5) [right] {$\sharp$};\node at(3.75,2.5) [right] {$\sharp$};
            \node at(0.5,0.25) [below] {$\sharp$};\node at(1.5,0.25) [below] {$\sharp$};\node at(2.5,0.25) [below] {$\sharp$};\node at(3.5,0.25) [below] {$\sharp$};
            \node at(0.5,2.75) [above] {$\sharp$};\node at(1.5,2.75) [above] {$\sharp$};\node at(2.5,2.75) [above] {$\sharp$};\node at(3.5,2.75) [above] {$\sharp$};
        }}
        $\rightarrow$
        \raisebox{-1.2cm}{\tikzpic{0.6}{
            \draw (0,0) grid (4,3);
            \node at(0.25,0.5) [left] {$\sharp$};\node at(0.25,1.5) [left] {$\sharp$};\node at(0.25,2.5) [left] {$\sharp$};
            \node at(3.75,0.5) [right] {$\sharp$};\node at(3.75,1.5) [right] {$\sharp$};\node at(3.75,2.5) [right] {$\sharp$};
            \node at(0.5,0.25) [below] {$\sharp$};\node at(1.5,0.25) [below] {$\sharp$};\node at(2.5,0.25) [below] {$\sharp$};\node at(3.5,0.25) [below] {$\sharp$};
            \node at(0.5,2.75) [above] {$\sharp$};\node at(1.5,2.75) [above] {$\sharp$};\node at(2.5,2.75) [above] {$\sharp$};\node at(3.5,2.75) [above] {$\sharp$};
            \node at(0.25,2) [right] {$\sharp$};\node at(1.25,2) [right] {$\sharp$};\node at(2.25,2) [right] {$\sharp$};\node at(3.25,2) [right] {$\sharp$};
            \node at(0.25,1) [right] {$1$};\node at(1.25,1) [right] {$1$};\node at(2.25,1) [right] {$\sharp$};\node at(3.25,1) [right] {$\sharp$};
            \node at(0.75,2.5) [right] {$1$};\node at(1.75,2.5) [right] {$\sharp$};\node at(2.75,2.5) [right] {$1$};
            \node at(0.75,1.5) [right] {$\sharp$};\node at(1.75,1.5) [right] {$\sharp$};\node at(2.75,1.5) [right] {$1$};
            \node at(0.75,0.5) [right] {$\sharp$};\node at(1.75,0.5) [right] {$\sharp$};\node at(2.75,0.5) [right] {$1$};
        }}.
    \end{center}

    When the tiling region is the entire plane (or an entire quadrant), the problem of whether a tiling exists to cover the entire plane (i.e., the tiling existence problem) is undecidable~\cite{WT, WT_1}.
    In fact, a significant portion of problems regarding Wang tilings on the plane are undecidable~\cite{WT_und_01}.

    When the tiling region is finite, the problem becomes decidable (since at least all possible tilings can be enumerated), but it remains computationally complex.
    For example, the existence problem for Wang tilings on a $l \times h$ chessboard (whether or not boundary characters are restricted or tile rotation is allowed)
    is always NP, while the enumeration problem is always \#P~\cite{WT_CC, WT_CC_1}.
    Therefore, only approximate enumeration algorithms can be provided for large-scale general Wang tilings.

    $n$\textit{-ominoes} are the natural generalization of dominoes, referring to the tiles composed of $n$ unit squares joined edge-to-edge~\cite{nOT_1, TC}.

    $3$-ominoes are also called \textit{trominoes}, with $2$ shapes (as shown below):
    the \textit{right tromino} and the \textit{straight tromino} ~\cite{LOT}.

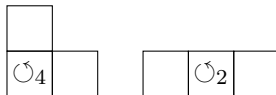
    \begin{center}
        \tikzpic{0.6}{
            \drawtromino{0}{0}{0}{1}{1}{0};
            \node at(0.5,0.5) {$\circlearrowleft_4$};
            \drawtromino{3}{0}{1}{0}{2}{0};
            \node at(4.5,0.5) {$\circlearrowleft_2$};
        }
    \end{center}
    The notation $\circlearrowleft_n$ on a tile indicates that the polyomino can be rotated into $n$ different orientations.

    Similarly, $4, 5, \ldots$-ominoes have special names: tetrominoes, pentominoes, $\ldots$.

    \cite{til} provides a method for simulating Wang tilings using polyomino tilings, proving that the problem of whether general domino tilings can cover the entire plane is undecidable,
    and that the existence and enumeration problems for tilings on general chessboards are NP and \#P, respectively.

    Currently, for problems of this type, backtracking algorithms or heuristic searches are commonly employed to exhaustively enumerate all possible tilings~\cite{DT_alg, nOT_enum, til_alg_2, til_alg_3}.
    However, due to the exponential growth of computational cost with increasing problem size, such computations may require supercomputers or distributed computing systems, and in some cases, may even become entirely infeasible.

    The algorithm presented in this paper is an optimized \textit{dynamic programming} approach, which does not require the exhaustive enumeration of all possible tilings.
    We store the states by using multi-dimensional tensors (referred to as \textit{state tensors}), and state transitions are implemented through specific tensor operations, which can significantly save space and improving efficiency.
    Furthermore, we employ modular enumeration and solving congruence equations instead of direct enumeration to avoid data overflow and further improve computational efficiency.

    \subsection{Tensor contractions}

    This paper uses $\bbK$ to denote a general commutative (semi-)ring.
    $a_i$ denotes a free index from the index set $\varSigma_i$ (also referred to as the \textit{vertical alphabet}), and $\sigma_i$ denotes its size;
    $\alpha_i$ denotes a free index from the index set $\varTheta_i$ (also referred to as the \textit{horizontal alphabet}), and $\theta_i$ denotes its size.

    This paper uses \textit{string diagrams} to represent tensors and their operations~\cite{TN, TN_1}.
    Specifically, we use closed shapes to represent tensors and pins to represent tensor indices.
    For example, a tensor $T_{a_1 a_2 a_3 a_4}$ with $4$ indices is represented as

    \begin{center}
        \raisebox{-1cm}{\tikzpic{0.6}{
            \drawtile{0}{0};\node at(0.5,0.5) {$T$};
            \draw (0,0.5)--(-0.5,0.5);\draw (1,0.5)--(1.5,0.5);\draw (0.5,0)--(0.5,-0.5);\draw (0.5,1)--(0.5,1.5);
            \node at(-0.5,0.5) [left] {$a_1$};
            \node at(0.5,1.5) [above] {$a_2$};
            \node at(0.5,-0.5) [below] {$a_3$};
            \node at(1.5,0.5) [right] {$a_4$};
            \node at(0,0) {$\circ$};\draw[->] (0,0)--(0,0.3);
        }}
        or
        \raisebox{-0.5cm}{\tikzpic{0.6}{
            \draw (0,0) rectangle (4,0.5);\node at(2,0.25) {$T$};
            \draw (0.5,0.5)--(0.5,1);\draw (1.5,0.5)--(1.5,1);\draw (2.5,0.5)--(2.5,1);\draw (3.5,0.5)--(3.5,1);
            \node at(0.5,1) [above] {$a_1$};\node at(1.5,1) [above] {$a_2$};\node at(2.5,1) [above] {$a_3$};\node at(3.5,1) [above] {$a_4$};
            \node at(0,0.5) {$\circ$};\draw[->] (0,0.5)--(0.3,0.5);
        }}.
    \end{center}

    The order of the indices (from left to right) is determined by starting at the hollow points on the graph and following the direction of the arrows.  
    According to this order, we can define \textit{the order of index sequences} as \textit{the lexicographically order from left to right}
    (i.e., from left to right as low-order to high-order).
    For example, when each alphabet is $\{0, 1, 2\}$,
    \begin{align*}
              & 000\ldots 0 \prec 100\ldots 0 \prec 200\ldots 0
        \prec 010\ldots 0 \prec 110\ldots 0 \prec 210\ldots 0 \\
        \prec~& 020\ldots 0 \prec 120\ldots 0 \prec 220\ldots 0
        \prec 001\ldots 0 \prec 101\ldots 0 \prec 201\ldots 0 \\
        \prec~& \ldots \prec 000\ldots 1 \prec \ldots \prec 000\ldots 2 \prec \ldots.
    \end{align*}
    When the $i$-th alphabet is taken as the set of natural numbers
    \[
        \varSigma_i = \{0, 1, \ldots, (\sigma_i - 1)\},
    \]
    the index of each index sequence is
    \[
        a_1 a_2 a_3 \ldots a_l \mapsto a_1+\sigma_1 a_2 +\sigma_1 \sigma_2 a_3+\ldots+\sigma_1 \sigma_2\ldots\sigma_{c-1}a_l.
    \]
    The general way computers store tensors is to store the values corresponding to the index sequences in this order.
    It is important to note that in programming languages, when accessing elements in multidimensional lists, the order is generally from left to right as \textit{high-order to low-order},
    resulting in index sequences that are exactly the reverse of what we define here.
    For example,
    \[
        T_{a_1 a_2 a_3 a_4} \mapsto \mathrm{T[a\_4][a\_3][a\_2][a\_1]}.
    \]
    The shape parameter of the tensor (containing the sizes of the index sets for each index) is also reversed in this way.

    Additionally, we define the rotation of indices as
    \[
        T_{a_1 a_2 \ldots a_l} \mapsto T_{a_2 \ldots a_l a_1},
    \]
    which means rotating the first index to the last position.
    If we specify the order of indices by assigning an initial point and an initial direction on the graphical representation of the tensor, such a rotation can be achieved by moving the initial point forward and jumping over an index, like

    \begin{center}
        $T_{a_1 a_2 \ldots a_l} \mapsto T_{a_2 \ldots a_l a_1} ~~\Leftrightarrow$
        \raisebox{-1cm}{\tikzpic{0.6}{
            \drawtile{0}{0};\node at(0.5,0.5) {$T$};
            \draw (0,0.5)--(-0.5,0.5);\draw (1,0.5)--(1.5,0.5);\draw (0.5,0)--(0.5,-0.5);\draw (0.5,1)--(0.5,1.5);
            \node at(-0.5,0.5) [left] {$a_1$};
            \node at(0.5,1.5) [above] {$a_2$};
            \node at(0.5,-0.5) [below] {$a_3$};
            \node at(1.5,0.5) [right] {$a_4$};
            \node at(0,0) {$\circ$};\draw[->] (0,0)--(0,0.3);
        }}
        $\mapsto\!\!$
        \raisebox{-1cm}{\tikzpic{0.6}{
            \drawtile{0}{0};\node at(0.5,0.5) {$T$};
            \draw (0,0.5)--(-0.5,0.5);\draw (1,0.5)--(1.5,0.5);\draw (0.5,0)--(0.5,-0.5);\draw (0.5,1)--(0.5,1.5);
            \node at(-0.5,0.5) [left] {$a_1$};
            \node at(0.5,1.5) [above] {$a_2$};
            \node at(0.5,-0.5) [below] {$a_3$};
            \node at(1.5,0.5) [right] {$a_4$};
            \node at(0,1) {$\circ$};\draw[->] (0,1)--(0.3,1);
        }}.
    \end{center}

    By dividing the indices into row indices (the subscripts of the tensor, represented by outward arrows in the string diagram) and column indices (the superscripts, represented by inward arrows),
    a tensor can be flattened into a matrix, like

    \begin{center}
        $T_{a_1 a_2 a_3 a_4}\mapsto T_{a_1 a_2}^{a_3 a_4} ~~\Leftrightarrow$
        \raisebox{-1cm}{\tikzpic{0.6}{
            \drawtile{0}{0};\node at(0.5,0.5) {$T$};
            \draw (0,0.5)--(-0.5,0.5);\draw (1,0.5)--(1.5,0.5);\draw (0.5,0)--(0.5,-0.5);\draw (0.5,1)--(0.5,1.5);
            \node at(-0.5,0.5) [left] {$a_1$};
            \node at(0.5,1.5) [above] {$a_2$};
            \node at(0.5,-0.5) [below] {$a_3$};
            \node at(1.5,0.5) [right] {$a_4$};
            \node at(0,0) {$\circ$};\draw[->] (0,0)--(0,0.3);
        }}
        $\mapsto\!\!$
        \raisebox{-1cm}{\tikzpic{0.6}{
            \drawtile{0}{0};\node at(0.5,0.5) {$T$};
            \draw[->] (0,0.5)--(-0.5,0.5);\draw[<-] (1,0.5)--(1.5,0.5);\draw[<-] (0.5,0)--(0.5,-0.5);\draw[->] (0.5,1)--(0.5,1.5);
            \node at(-0.5,0.5) [left] {$a_1$};
            \node at(0.5,1.5) [above] {$a_2$};
            \node at(0.5,-0.5) [below] {$a_3$};
            \node at(1.5,0.5) [right] {$a_4$};
            \node at(0,0) {$\circ$};\draw[->] (0,0)--(0,0.3);
        }}
        $=$
        \raisebox{-0.8cm}{\tikzpic{0.6}{
            \draw (0,0) rectangle (2,0.5);\node at(1,0.25) {$T$};
            \draw[->] (0.5,0.5)--(0.5,1);\draw[->] (1.5,0.5)--(1.5,1);\draw[->] (0.5,-0.5)--(0.5,0);\draw[->] (1.5,-0.5)--(1.5,0);
            \node at(0.5,1) [above] {$a_1$};\node at(1.5,1) [above] {$a_2$};
            \node at(0.5,-0.5) [below] {$a_2$};\node at(1.5,-0.5) [below] {$a_3$};
            \node at(0,0.5) {$\circ$};\draw[->] (0,0.5)--(0.3,0.5);
        }}.
    \end{center}

    For programming languages where the matrix storage format is C-style (row-major, such as C and Python), if the flattened row and column indices retain their original order,
    and the row indices precede the column indices (as in the example above),
    the flattening operation does not modify the data storing the tensor or matrix, only the shape parameters, so the time cost is negligible.
    Note that some languages use F-style matrix storage (column-major, such as Matlab).

    After flattening both tensors into matrices, the contraction of tensor indices can be represented as matrix multiplication, like

    \begin{center}
        $(T_1)_{a_1 a_2}^{a_3 a_4} (T_2)_{a_3 a_4}^{a_5 a_6}= T_{a_1 a_2}^{a_5 a_6} ~~\Leftrightarrow $
        \raisebox{-1.4cm}{\tikzpic{0.6}{
            \draw (0,0) rectangle (2,0.5);\node at(1,0.25) {$T_1$};
            \draw[->] (0.5,0.5)--(0.5,1);\draw[->] (1.5,0.5)--(1.5,1);\draw[->] (0.5,-1.5)--(0.5,0);\draw[->] (1.5,-1.5)--(1.5,0);
            \node at(0.5,1) [above] {$a_1$};\node at(1.5,1) [above] {$a_2$};
            \node at(0.5,-0.75) {$a_3$};\node at(1.5,-0.75) {$a_4$};
            \node at(0,0.5) {$\circ$};\draw[->] (0,0.5)--(0.4,0.5);

            \draw (0,-2) rectangle (2,-1.5);\node at(1,-1.75) {$T_2$};
            \draw[->] (0.5,-2.5)--(0.5,-2);\draw[->] (1.5,-2.5)--(1.5,-2);
            \node at(0.5,-2.5) [below] {$a_5$};\node at(1.5,-2.5) [below] {$a_6$};
            \node at(0,-1.5) {$\circ$};\draw[->] (0,-1.5)--(0.4,-1.5);
        }}
        $~~=$
        \raisebox{-0.8cm}{\tikzpic{0.6}{
            \draw (0,0) rectangle (2,0.5);\node at(1,0.25) {$T$};
            \draw[->] (0.5,0.5)--(0.5,1);\draw[->] (1.5,0.5)--(1.5,1);\draw[->] (0.5,-0.5)--(0.5,0);\draw[->] (1.5,-0.5)--(1.5,0);
            \node at(0.5,1) [above] {$a_1$};\node at(1.5,1) [above] {$a_2$};
            \node at(0.5,-0.5) [below] {$a_5$};\node at(1.5,-0.5) [below] {$a_6$};
            \node at(0,0.5) {$\circ$};\draw[->] (0,0.5)--(0.3,0.5);
        }},
    \end{center}
    where $a_3,a_4$ are the contracted indices.

    \section{Wang tiling}

    \subsection{In helicoids}

    We first discuss the enumeration (or weighted enumeration) problem of Wang tilings in a \textit{helicoid} (which is called a \textit{twisted cylinder} in~\cite{TC}),
    as its properties are more convenient for demonstrating the principles of our algorithm and are sufficiently general.

    First, we generalize the classical Wang tile set to a $\bbK$-Wang tile set, where each tile carries a value from $\bbK$ as its \textit{weight},
    and the weights are multiplied during tiling.
    In other words, we extend the set of Wang tiles to the $bbK$-module generated by Wang tiles.
    For example, the horizontal concatenation of two $\bbK$-Wang tiles is
    \[
        (k_1\cdot \wtile{$\alpha_1$}{$a_1$}{$b_1$}{$\beta_1$})\circ (k_2\cdot\wtile{$\alpha_2$}{$a_2$}{$b_2$}{$\beta_2$})
        =\begin{cases}
             k_1 k_2\cdot\whdomino{$\alpha_1$}{$a_1$}{$b_1$}{$\beta_1$}{}{$a_2$}{$b_2$}{$\beta_2$}, &\beta_1=\alpha_2;\\
             0, & \beta_1\neq \alpha_2,
        \end{cases}
    \]
    where $0$ is the zero element of $\bbK$.
    The total weight after tiling is defined as the sum of all possible tile pairings.
    Thus, the original Wang tile set and its operations can be seen as a special case when $\bbK$ is taken as the Boolean semiring $\left< \{0,1\}, \vee, \wedge \right>$.

    It can be observed that there is an obvious correspondence between $\bbK$-Wang tile sets and $\bbK$-tensors with $4$ indices,
    where the characters on the four edges are treated as indices, and the weight of the tile is treated as the value in the tensor corresponding to those indices, like
    \[
        k\cdot \wtile{$\alpha$}{$a$}{$b$}{$\beta$} \mapsto T_{b \beta\alpha a}=k
    \]
    Therefore, both horizontal and vertical tiling of $\bbK$-Wang tiles can be seen as tensor index contractions, and the rules are completely equivalent.

    When $\bbK$ is taken as $\bbN$ or $\bbZ$, the weight of a tile can be understood as the \textit{enumeration},
    and a tile with a weight $k$ means that there are $k$ tiles in the tile set with the same edge characters but distinguishable from each other (e.g., with different central characters or colors).
    In this case, calculating the total weight after the tensor contraction under a given order is equivalent to enumerating the tilings under this order.

    The construction of the helicoid discussed in this section is as follows: take a strip composed of $S$ squares connected horizontally, and wrap it around a cylinder with circumference $l$,
    such that the square with index $i + l$ is placed below the square with index $i$ (as shown below, the numbers in the squares represent their indices).

    \begin{center}
        \tikzpic{0.5}{
            \draw (1,0) grid (3,4);\draw (3,1) grid (6,5);
            \draw (1,4.5) arc (90:270:0.25);\draw (0.75,0.25) arc (180:270:0.25);
            \draw (6,5) arc (90:0:0.25);\draw (6,1) arc (90:-90:0.25);\draw (3,0.5)--(6,0.5);
            \draw (0.75,0.25)--(0.75,4.25);\draw (6.25,0.75)--(6.25,4.75);\draw (1,4.5)--(3,4.5);
            \node at(3.5,4.5) {$1$};\node at(4.5,4.5) {$2$};\node at(5.5,4.5) {$3$};\node at(1.5,3.5) {$\ldots$};
            \node at(2.5,3.5) {$l$};\node at(3.5,3.5) {\tiny$l\!+\!1$};\node at(4.5,3.5) {\tiny$l\!+\!2$};\node at(5.5,3.5) {$\ldots$};
            \node at(1.5,0.5) {$\ldots$};\node at(2.5,0.5) {$S$};
        }
        ~~~~~~~~
        \tikzpic{0.5}{
            \draw (1,0) grid (3,1);\draw (5,1) grid (6,2);
            \draw (1,1.5) arc (90:270:0.25);\draw (0.75,0.25) arc (180:270:0.25);
            \draw (6,2) arc (90:0:0.25);\draw (6,1) arc (90:-90:0.25);\draw (3,0.5)--(6,0.5);
            \draw (0.75,0.25)--(0.75,1.25);\draw (6.25,0.75)--(6.25,1.75);\draw (1,1.5)--(5,1.5);
            \draw[line width=2pt] (3,1)--(5,1);
        }
    \end{center}

    It is not difficult to see that when $S \geq l$ (as shown in the left), the top and bottom edges of the helicoid always have the same shape: a twisted line of length $l$, connected at both ends by a unit vertical line.
    For the degenerate case where $S < l$ (as shown in the right), the boundary of the strip can be extended (i.e., adding the bold lines in the figure) to maintain this edge shape.

    If tiles are placed in the helicoid in the order of the original horizontal strip, the remaining region will still be a helicoid with circumference $l$,
    and the shape of the edges will not change.

    We assign a character from the vertical alphabet $\varSigma$ to each horizontal edge of the helicoid's top edge (ignoring the bottom edge for now), and a character from the horizontal alphabet $\varTheta$ to each vertical edge,
    and always treat \textit{the lower endpoint of the vertical edge as the starting point}.
    The characters on the vertical edge (as indices) are arranged \textit{counterclockwise} to form a word (sequence) of length $l + 1$, taken from $\varTheta\varSigma^l$.
    The combination of such words with weights from $\bbK$ ($\bbK$-word set) can be represented as tensors in the same way as $\bbK$-tile sets, like

    \begin{center}
        \raisebox{-0.5cm}{\tikzpic{0.5}{
            \draw (1,2) arc (90:270:1 and 1);
            \draw (6,1) arc (-90:90:0.5 and 0.5);
            \draw (1,0)--(1,1)--(6,1);\draw (1,2)--(6,2);
            \node at(0.5,1.5) {$W$};
            \node at(1,0.5) {$\alpha_1$};\node at(1.5,1) {$a_1$};
            \node at(2.5,1) {$a_2$};\node at(3.5,1) {$\ldots$};\node at(5.5,1) {$a_l$};
        }}
        $\mapsto$
        \raisebox{-0.5cm}{\tikzpic{0.5}{
            \draw (0,-0.5)--(0,1)--(6,1)--(5.5,1.5)--(-1,1.5)--(-1,0)--(0,-0.5);
            \node at(-0.5,0.5) {$W$};
            \draw (0,0)--(0.5,0);\draw (1.5,1)--(1.5,0.5);\draw (2.5,1)--(2.5,0.5);\draw (3.5,1)--(3.5,0.5);\draw (5.5,1)--(5.5,0.5);
            \node at(0.5,0) {$\alpha_1$};\node at(1.5,0.5) [below] {$a_1$};
            \node at(2.5,0.5) [below] {$a_2$};\node at(3.5,0.5) [below] {$\ldots$};\node at(5.5,0.5) [below] {$a_l$};
            \node at(0,-0.5) {$\circ$};
            \draw[->] (0,-0.5)--(0,-0.1);
        }}.
    \end{center}

    Next, we take $\bbK$ as $\bbZ$ and discuss the enumeration algorithm for such tilings.
    Consider the $\bbK$-word set (tensor representation) of the word on the top edge of the helicoid at a certain moment as $W_{\alpha_1 a_1 a_2 \ldots a_l}$.
    At this time, place a Wang tile $T_{b_1\alpha_2\alpha_1 a_1}$ at the first position of the helicoid region (that is, the $\Gamma$-shaped bend on the top edge).
    The enumeration of the new word on the edge (note that the starting point has moved one unit to the right) is
    \[
        (W\ltimes T)_{\beta a_2\ldots a_l b_1}  \coloneq \sum_{\alpha a_1\in \varTheta\varSigma}
        W_{\alpha a_1 a_2 \ldots a_l} T_{b_1\alpha_2\alpha_1 a_1},
    \]
    and the string diagram is
    \begin{center}
        \raisebox{-1.5cm}{\tikzpic{0.5}{
            \draw (0,-0.5)--(0,1)--(5,1)--(4.5,2)--(-1,2)--(-1,0)--(0,-0.5);
            \node at(-0.5,0.5) {$W$};
            \draw (0,0)--(0.5,0);\draw (1.5,1)--(1.5,0.5);
            \draw (2.5,1)--(2.5,0.5);\draw (3.5,1)--(3.5,0.5);\draw (4.5,1)--(4.5,0.5);
            \node at(2.5,0.5) [below] {$a_2$};\node at(3.5,0.5) [below] {$\ldots$};\node at(4.5,0.5) [below] {$a_l$};
            \node at(0,-0.5) {$\circ$};
            \draw[->] (0,-0.5)--(0,-0.1);

            \draw (1,-1.5) rectangle (2,-0.5);
            \node at(1.5,-1) {$T$};
            \draw (1,-1.0)--(0.5,-1.0)--(0.5,0);\draw (1.5,-0.5)--(1.5,0.5);
            \node at(0.5,0) {$\alpha_1$};\node at(1.5,0.25) {$a_1$};
            \draw (2.5,-1.0)--(2,-1.0);\draw (1.5,-2)--(1.5,-1.5);
            \node at(2.5,-1.0) [right] {$\alpha_2$};\node at(1.5,-2) [below] {$b_1$};
        }}
        $\mapsto$
        \raisebox{-1cm}{\tikzpic{0.5}{
            \draw (0,-0.5)--(1,-0.5)--(1,1)--(5,1)--(4.5,2)--(-1,2)--(-1,0)--(0,-0.5);
            \node at(0,1.5) {$W\!\ltimes \! T$};
            \draw (1,0)--(1.5,0);\draw (0.5,-0.5)--(0.5,-1);
            \node at(1.5,0) {$\alpha_2$};\node at(0.5,-1) [below] {$b_1$};
            \draw (2.5,1)--(2.5,0.5);\draw (3.5,1)--(3.5,0.5);\draw (4.5,1)--(4.5,0.5);
            \node at(2.5,0.5) [below] {$a_2$};\node at(3.5,0.5) [below] {$\ldots$};\node at(4.5,0.5) [below] {$a_l$};
            \node at(1,-0.5) {$\circ$};
            \draw[->] (1,-0.5)--(1,-0.1);
            \draw[->] (0,-0.3) arc (120:75: 1);
        }}
        $=$
        \raisebox{-0.5cm}{\tikzpic{0.5}{
            \draw (0,-0.5)--(0,1)--(5,1)--(4.5,2)--(-1,2)--(-1,0)--(0,-0.5);
            \node at(0,1.5) {$W\!\ltimes \! T$};
            \draw (0,0)--(0.5,0);\draw (1.5,1)--(1.5,0.5);\draw (2.5,1)--(2.5,0.5);\draw (3.5,1)--(3.5,0.5);\draw (4.5,1)--(4.5,0.5);
            \node at(0.5,0) {$\alpha_2$};\node at(1.5,0.5) [below] {$a_2$};
            \node at(2.5,0.5) [below] {$\ldots$};\node at(3.5,0.5) [below] {$a_l$};\node at(4.5,0.5) [below] {$b_1$};
            \node at(0,-0.5) {$\circ$};
            \draw[->] (0,-0.5)--(0,-0.1);
        }}.
    \end{center}
    That is, the indices $\alpha_1 a_1$ on the seam are contracted, and then the index $b_1$ corresponding to the character below the tile is rotated to the end (corresponding to the movement of the starting point).
    Without changing the encoding, the above can also be expressed in the form of matrix operations, such as
    \[
        (W\ltimes T)_{\alpha_2 a_2}^{a_3\ldots a_l b_1}  \coloneq
        W_{\alpha_1 a_1}^ {a_2\ldots a_l} T_{ b_1\alpha_2}^{\alpha_1 a_1}.
    \]
    However, directly multiplying the two matrices on the right would result in a matrix with indices $(WT)_{b_1 \alpha_2}^{a_2 \ldots a_l}$,
    which differs from our defined $W \ltimes T$ by an index rotation.
    The reason for defining $W \ltimes T$ in this way is to ensure that in each matrix multiplication, the column indices (the subscripts) of the left matrix are always the first two indices,
    thus ensuring that the index contractions can be converted into matrix multiplications.

    To save the computational cost of rotating indices, we transform the above operations into equivalent batched matrix multiplications.

    Note that if the free index $b_1$ is fixed to $b^* \in \varSigma$ (i.e., taking a slice of the tensor along this index), then
    \[
        (W\ltimes T)_{\alpha_2 a_2\ldots a_l b^*}  =
        W_{\alpha_1 a_1}^ {a_2\ldots a_l} T_{ b^*\alpha_2}^{\alpha_1 a_1},
    \]
    where the encoding of $T_{b^* \alpha_2}^{\alpha_1 a_1}$ equals to a slice of the tensor $T_{b_1 \alpha_2\alpha_1 a_1}$.
    Let $\varSigma=\{0,1,\ldots,(\sigma-1)\}$, according to the tensor encoding rules, the encoding of $(W \ltimes T)_{\alpha_2 a_2 \ldots a_l b_1}$ equals to
    the concatenated encoding of
    \[
        (W\ltimes T)_{\alpha_2 a_2\ldots a_l 0}, (W\ltimes T)_{\alpha_2 a_2\ldots a_l 1},
        \ldots, (W\ltimes T)_{\alpha_2 a_2\ldots a_l (\sigma-1)}.
    \]
    Therefore, to obtain the encoding of $(W \ltimes T)$, we only need to multiply $W_{\alpha_1 a_1}^{a_2 \ldots a_l}$ sequentially with
    \[
        T_{0\alpha_2}^{\alpha_1 a_1}, T_{1\alpha_2}^{\alpha_1 a_1},\ldots,T_{(\sigma-1)\alpha_2}^{\alpha_1 a_1},
    \]
    and concatenate the resulting matrices vertically (this concatenation does not consume time, as it only requires connecting the output memory spaces).

    When using NumPy, we can represent $T$ as a $3$-dimensional tensor with shape $(\sigma, \sigma\theta, \theta)$,
    which is obtained by stacking $\sigma$ (corresponding to the first value of the shape parameter) matrices of size $\sigma\theta \times \theta$ (corresponding to the last two values)
    \[
        T_{0 \alpha_2}^{\alpha_1 a_1}, T_{1 \alpha_2}^{\alpha_1 a_1}, \ldots, T_{(\sigma-1) \alpha_2}^{\alpha_1 a_1},
    \]
    along the third dimension.
    Thus, when computing $W \ltimes T$, we only need to reshape $W$ to $(\sigma^{l-1}, \sigma\theta)$ and perform NumPy's matrix multiplication.
    According to its broadcasting mechanism, since $W$ lacks one dimension, it will first be expanded by adding a dimension and broadcasted (copied into $\sigma$ copies and stacked),
    resulting in a $3$-dimensional tensor with shape $(\sigma, \sigma^{l-1}, \sigma\theta)$.
    Then, the last two dimensions will be treated as the number of rows and columns for matrix multiplication, yielding a $3$-dimensional tensor with shape $(\sigma, \sigma^{l-1}, \theta)$,
    which is exactly the tensor obtained by stacking the matrices
    \[
        (W \ltimes T)_{\alpha_2}^{a_2 \ldots a_l 0}, (W \ltimes T)_{\alpha_2}^{a_2 \ldots a_l 1},
        \ldots, (W \ltimes T)_{\alpha_2}^{a_2 \ldots a_l (\sigma-1)},
    \]
    along the third dimension.
    For languages without such broadcasting mechanisms, the above operations are also straightforward to implement.

    Starting from the initial top-edge word set that only contain $\sharp\sharp\ldots\sharp$ (with a weight $1$),
    using this algorithm to place tiles sequentially and summing the results, we can
    calculate the number $N$ of the word on the bottom edge after tiling the entire helicoid.
    Different positions can also specify different tile sets, i.e., defining
    \[
        W_{i+1} = (W_i \ltimes T_i),
    \]
    where $T_i$ is the $\bbZ$ tile set at the position with index $i$.

    However, when the word enumeration becomes too large, matrix operations are prone to data overflow.
    To address this issue, we select multiple prime numbers $p$ (discussed later),
    and modify the above operations to modulo $p$ operations, i.e., taking the result modulo $p$ after each matrix operation.
    Thus, the tiling enumeration becomes a modulo $p$ enumeration, projecting the ring $\bbZ$ of the tensor onto $\bbZ_p$, and the $\bbZ$-tensor is projected onto a $\bbZ_p$-tensor.
    It is easy to prove that such a projection preserves all the tensor operations used above.

    According to the Chinese Remainder Theorem, if $\gcd(p, P) = 1$, then the system of congruences for $R'$ (where $[R']_p$ denotes the congruence class of $R'$ modulo $p$)
    \[
        \begin{cases}
            [R']_p = [r]_p; \\
            [R']_P = [R]_P
        \end{cases}
    \]
    has the solution
    \[
        [R']_{pP} = [P^- P r + p^- p R]_{pP},
    \]
    where $P^- \in [P]_p^{-1}$ and $p^- \in [p]_P^{-1}$ satisfy
    \[
        P^- P + p^- p = 1.
    \]
    Such $P^-$ and $p^-$ can be obtained using the Euclidean algorithm on $(p, P)$.

    Thus, the following steps can be used to compute the enumeration in the general sense through modular enumeration:
    First, estimate the lower and upper bounds of the desired enumeration, such as
    \[
        \underline{N} \leq N < \overline{N}.
    \]
    Then, assume $N$ lies in the congruence class $[R]_p$.
    Starting with $R = 0$ and $P = 1$ (where $[0]_1 = \bbZ$), select sufficiently many prime numbers (this ensures $\gcd(p, P) = 1$) $p$ for modular enumeration,
    and replace the original $[R]_P$ with the obtained $[R']_{pP}$, repeating the process as follows:

    \[
        \begin{cases}
            P \mapsto pP; \\
            R \mapsto (P^- P r + p^- pR) ~\%~ pP.
        \end{cases}
    \]

    Here, $\%$ denotes the modulo operation, ensuring $R < P$.
    After repeating this process a sufficient number of times until $P > \overline{N} - \underline{N}$, there will be only one number in $[R]_P$ that falls within $[\underline{N}, \overline{N})$,
    which is the desired enumeration.
    This number can be obtained using
    \[
        N = \underline{N} + (R - \underline{N}) ~\%~ P.
    \]

    When the product $P$ of the accumulated primes is less than $(\overline{N} - \underline{N})$, the above formula may not yield the desired enumeration;
    however, when $P \geq (\overline{N} - \underline{N})$, $N$ will be the required enumeration and will remain fixed, with no further changes to $R$ or $P$ affecting its value.

    In this case, to make the algorithm stop faster, a \textit{halting mechanism} can be added: given an integer $\tau > 0$,
    stop when the above formula yields the same $N$ for $\tau$ consecutive rounds of solving the congruence equations, and output this $N$ as the result.
    Note that $R$ does not change if and only if
    \[
        p\mid (P^-Pr + p^-pR) \Leftrightarrow p\mid r \Leftrightarrow p\mid N,
    \]
    i.e., $p$ exactly divides the desired enumeration, with a probability (though not strictly mathematical) of approximately $\frac{1}{p}$.

    Therefore, under the above setup, the probability of misjudgment is approximately the reciprocal of the product of the last $\tau$ primes used.
    At the same time, if a misjudgment occurs, the deviation is at least the product of these $\tau$ primes.
    As long as the probability of misjudgment is acceptable, or the deviation is large enough to be detected (if an error is detected, increase $\tau$ and recompute),
    the above setup can be used to stop the computation more promptly.

    Next, we discuss the selection of prime numbers $p$.
    If $p$ is chosen too large, it may still cause data overflow, but under the premise of avoiding overflow, the larger $p$ is, the more efficient the algorithm becomes.
    Therefore, we need to determine the upper bound $\overline{p}$ for the selection of $p$, and search for primes starting from this number downward to use as $p$.

    Note that each matrix multiplication in the algorithm takes the form
    \[
        W_{\alpha_1 a}^{\ldots} T_{b^* \alpha_2}^{\alpha_1 a},
    \]
    i.e., the rows of $W_{\alpha_1 a}^{\ldots}$ are multiplied by the columns of $T_{b^* \alpha_2}^{\alpha_1 a}$.
    Here, the entries of $W_{\alpha_1 a}^{\ldots}$ are uncertain, but after taking modulo $p$, they do not exceed $p - 1$.
    There are two cases to consider:

    (1) The terms of $T$ are all non-negative.

    In this case, assume the data type is a 32-bit unsigned integer, with a range from $0$ to $2^{32} - 1$.
    In a tiling enumeration where the tile sets at each position are $T_i$, we define
    \[
        S_{\mathrm{col}} \coloneq \max_{i} \max_{\alpha_2} \sum_{a \alpha_1} {(T_i)}_{b^* \alpha_2}^{\alpha_1 a},
    \]
    as the maximum column sum of $(T_i)_{0 \alpha_2}^{\alpha_1 a_1}$ during matrix multiplication.
    Then, the largest possible number obtained from row-column multiplication should satisfy
    \[
        (p - 1) S_{\mathrm{col}} \leq 2^{32} - 1.
    \]
    Thus, $p$ only needs to be no larger than the upper bound
    \[
        \overline{p} = \left\lfloor \frac{2^{32} - 1}{S_{\mathrm{col}}} \right\rfloor + 1
    \]
    to ensure no data overflow occurs.

    (2) $T$ has negative terms.

    In this case, assume the data type is a 32-bit signed integer, with a range from $-2^{31}$ to $2^{31} - 1$.
    We need to ensure that the data does not overflow in either the positive or negative direction.

    Let $T_i^+$ denote the tensor obtained by retaining only the positive entries of $T_i$ (setting all other entries to $0$),
    and let $T_i^-$ denote the tensor obtained by retaining only the negative entries.
    Define
    \begin{align*}
        S^+_{\mathrm{col}} & \coloneq \max_{i} \max_{\alpha_2} \sum_{a \alpha_1} {(T^+_i)}_{b^* \alpha_2}^{\alpha_1 a}; \\
        S^-_{\mathrm{col}} & \coloneq \max_{i} \max_{\alpha_2} \sum_{a \alpha_1} {(T^-_i)}_{b^* \alpha_2}^{\alpha_1 a}.
    \end{align*}
    Then, the largest and smallest possible numbers obtained from row-column multiplication should satisfy
    \begin{align*}
        (p - 1) S^+_{\mathrm{col}} & \leq 2^{31} - 1; \\
        (p - 1) S^-_{\mathrm{col}} & \geq -2^{31}.
    \end{align*}
    Thus, the upper bound for $p$ should be
    \[
        \overline{p} = \left\lfloor \min\left(\frac{2^{31} - 1}{S^+_{\mathrm{col}}}, -\frac{2^{31}}{S^-_{\mathrm{col}}}\right) \right\rfloor + 1.
    \]

    \subsection{In rectangle chessboards}

    In this section, we consider the following Wang tiling in a chessboard:
    Tiling a $l \times h$ (width and height, respectively) chessboard with Wang tiles, where each position on the boundary of the region has the character $\sharp$ (no longer shown in the figures).
    Note that a rectangle can always be rolled into a helicoid by
    horizontally joining the pair of positions with horizontal and vertical indices $(l, j)$ and $(1, j+1)$ together, like

    \begin{center}
        \raisebox{-0.9cm}{\tikzpic{0.6}{
            \draw (0,1) grid (5,5);
            \node at(0.5,4.5) {$11$};\node at(1.5,4.5) {$21$};\node at(2.5,4.5) {$\ldots$};
            \node at(4.5,4.5) {$l1$};\node at(0.5,3.5) {$12$};\node at(1.5,3.5) {$22$};\node at(4.5,3.5) {$l2$};
            \node at(2.5,3.5) {$\ldots$};\node at(3.5,1.5) {$\ldots$};\node at(4.5,1.5) {$lh$};
        }}
        $~~\rightarrow$
        \raisebox{-1.5cm}{\tikzpic{0.6}{
            \draw[fill=gray!20] (1,1) rectangle (2,4);\draw[fill=gray!20] (1,1) rectangle (6,2);\draw[fill=gray!20] (5,1) rectangle (6,5);
            \draw[line width=2pt] (1,4)--(1,1)--(6,1)--(6,5);
            \draw (1,1) grid (6,5);
            \draw (1,4.5) arc (90:270:0.25);\draw (0.75,0.25) arc (180:270:0.25);
            \draw (6,5) arc (90:0:0.25);\draw (6,1) arc (90:-90:0.25);
            \draw (0.75,0.25)--(0.75,4.25);\draw (1,0)--(1,1);\draw (1,0.5)--(6,0.5);\draw (6.25,0.75)--(6.25,4.75);
            \node at(1.5,4.5) {$11$};\node at(2.5,4.5) {$21$};\node at(3.5,4.5) {$\ldots$};\node at(5.5,4.5) {$l1$};
            \node at(1.5,3.5) {$12$};\node at(2.5,3.5) {$22$};\node at(3.5,3.5) {$\ldots$};\node at(5.5,3.5) {$l2$};
            \node at(4.5,1.5) {$\ldots$};\node at(5.5,1.5) {$lh$};
        }}
        $~~=$
        \raisebox{-1.5cm}{\tikzpic{0.6}{
            \draw[fill=gray!20] (2,0) rectangle (3,4);\draw[fill=gray!20] (3,1) rectangle (6,2);\draw[fill=gray!20] (1,0) rectangle (3,1);\draw[fill=gray!20] (3,1) rectangle (4,4);
            \draw[line width=2pt] (3,4)--(3,1);\draw[line width=2pt] (1,0)--(3,0)--(3,1)--(6,1);
            \draw (1,0) grid (3,4);\draw (3,1) grid (6,5);
            \draw (1,4.5) arc (90:270:0.25);\draw (0.75,0.25) arc (180:270:0.25);
            \draw (6,5) arc (90:0:0.25);\draw (6,1) arc (90:-90:0.25);
            \draw (0.75,0.25)--(0.75,4.25);\draw (1,0)--(1,1);\draw (3,0.5)--(6,0.5);\draw (6.25,0.75)--(6.25,4.75);\draw (1,4.5)--(3,4.5);
            \node at(3.5,4.5) {$11$};\node at(4.5,4.5) {$21$};\node at(5.5,4.5) {$\ldots$};
            \node at(2.5,3.5) {$l1$};\node at(2.5,2.5) {$l2$};\node at(3.5,3.5) {$12$};\node at(4.5,3.5) {$22$};\node at(5.5,3.5) {$\ldots$};
            \node at(1.5,0.5) {$\ldots$};\node at(2.5,0.5) {$lh$};
        }}.
    \end{center}

    The Wang tiling in the chessboard can be viewed as Wang tiling in a helicoid with the initial top edge marked by $\sharp$,
    with the additional restriction that the characters on the bold lines in the figure must be $\sharp$.
    This can be achieved by modifying the tile sets at the boundary positions (shaded regions in the figure).
    For example, for all tile sets on the right edge, only tiles with the right edge marked by $\sharp$ are retained.
    Apply this processing operation separately to the left, right, and bottom edges.

    If we fix $l$ and want to quickly provide the enumeration for each height $h$, we can apply the above processing only to the left and right edges,
    and then perform tiling layer by layer.
    After each layer is fully tiled,
    note that the weight of the word $\sharp\sharp\ldots\sharp$ in the resulting $\bbZ_p$-word set is exactly a specific term of the state tensor at that moment.
    If $\sharp$ is the first character in the alphabet, then according to the encoding order, its weight is always the first term of the tensor.
    We only need to extract this term as the modulo enumeration of tilings satisfying the bottom edge condition at the current height.

    It is worth noting that since our algorithm does not require a common tile set for all positions, it also does not require the character sets at each position to be the same,
    only that the character sets on adjacent edges match.
    Therefore, after processing the tile sets on the edges as described above,
    we can retain the original character set on the edge or directly change the character set on that edge to $\{\sharp\}$.
    However, for simplicity, we usually fix the horizontal and vertical character sets to sizes $\theta$ and $\sigma$.

    Next, we discuss the complexity of the algorithm.
    At each moment, the $\bbZ_p$-word set of the twisted ring is taken from $\varTheta\varSigma^l$,
    so the state tensor (representing this $\bbZ_p$-word set) has a size of $\sigma^l \theta$.
    The algorithm needs to store two such tensors, occupying space
    $2s_\mathrm{data} \sigma^l \theta$
    where $s_\mathrm{data}$ is the memory size of a single integer variable (e.g., $s_\mathrm{data} = 4\mathrm{B}$ when the variable type is 32-bit integer).
    It can be seen that the space occupied by the state tensor increases extremely rapidly with $\sigma$ and $l$.
    When the data type of the integer is fixed, the total space complexity is
    \[
        \Theta(\sigma^l \theta).
    \]
    In comparison, the space occupied by other structures in the algorithm is negligible.

    Assuming $n_p$ is the number of primes need to be used, each prime requires tiling at $S = lh$ positions, and each tiling operation involves matrix multiplication,
    which requires traversing the state tensor once, resulting in a time complexity of
    \[
        \Theta(\sigma^l \theta n_p lh).
    \]

    If the desired enumeration is $N$, the distance between the upper and lower bounds is sufficiently large, and the automatic halting mechanism mentioned in the previous section is used to determine the number of primes selected,
    the algorithm will stop shortly after the product of the used primes exceeds $N$.
    The number of primes used can be estimated as
    \[
        n_p \approx \frac{\ln N}{\ln \overline{p}},
    \]
    where the constant $\overline{p}$ is determined by the selected $\bbZ$-tile set.
    And for a fixed $\bbZ$-tile set, the maximum value of $N$ obtained usually tends to increase exponentially with the increase of height $h$.
    Thus, we can also estimate the time complexity as
    \[
        \tilde{\Theta}(\sigma^l \theta lh^2).
    \]

    It can be seen that both the space and time complexity of the algorithm increase primarily with the size of the vertical character set $\sigma$ and the width of the chessboard $l$.
    The impact of other parameters on complexity is relatively small.
    When $l$ cannot be reduced, the size of the vertical character set should be minimized as much as possible to improve algorithm efficiency.
    Additionally, when $\sigma^l \theta > \theta^h \sigma$, to improve efficiency, the entire tiling region and tile sets should be transposed before enumerating.

    \subsection{In regions of other shapes}

    Our algorithm can also be further applied to the enumeration of Wang tilings in finite regions of other shapes.
    The approach remains to reduce the problem to Wang tilings on a helicoid.
    We will provide three examples.

    \begin{example}
        In irregular chessboards.
    \end{example}

    For an irregular chessboard, it can be viewed as a helicoid with certain blocks removed, which is equivalent to that each of these positions is occupied by a Wang tile \wtile{}{}{}{}.
    For example, consider the \textit{$3$rd order $3$-pillow} (\cite{til_enum_1}) shown below: restrict the characters on bold lines to be $\sharp$, and assign \wtile{}{}{}{} as the ($\bbK$-) Wang tile set to each of the gray positions,
    so the tiling enumeration in this area is converted into the tiling enumeration in a helicoid.

    \begin{center}
        \raisebox{-2.5cm}{\tikzpic{0.5}{
            \draw (0,6) grid (1,8);\draw (1,3) grid (2,9);\draw (2,0) grid (4,10);\draw (4,1) grid (5,7);\draw (5,2) grid (6,4);
        }}
        $~~\rightarrow$
        \raisebox{-2.5cm}{\tikzpic{0.5}{
            \draw (4,7) grid (5,9);\draw (5,4) grid (6,10);\draw (2,0) grid (4,10);\draw (4,1) grid (5,7);\draw (5,2) grid (6,4);
            \draw (6,10) arc (90:0:0.25);\draw (6,1) arc (90:-90:0.25);\draw (6.25,9.75)--(6.25,0.75);\draw (4,0.5)--(6,0.5);
            \draw (2,9.5) arc (90:270:0.25);\draw (1.75,0.25) arc (180:270:0.25);\draw (1.75,0.25)--(1.75,9.25);
            \draw[fill=gray] (5,1) rectangle (6,2);\draw[fill=gray] (4,9) rectangle (5,10);
            \draw[line width=2pt] (4,9)--(4,7)--(5,7)--(5,4)--(6,4)--(6,2);
            \draw[line width=2pt] (2,3)--(2,1);
            \draw[dashed] (2,3)--(1,3)--(1,6)--(0,6)--(0,8)--(1,8)--(1,9)--(2,9);
        }}
    \end{center}

    In fact, the enumeration of Wang tilings in irregular regions is equivalent to the contraction of tensors with irregular connectivity.
    Such calculation can often be optimized by adjusting the order of the tensor contraction, potentially leading to higher efficiency than our algorithm.
    For instance, consider the snake-shaped region shown below: if embedded in a helicoid, the size of the state tensor reaches a maximum of $\theta\sigma^{5}$.
    However, if the tensors are contracted in the order indicated by the arrows, the maximum size is only $\max(\theta, \sigma)$.
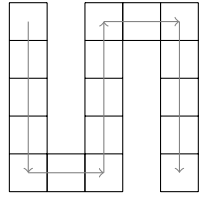
    \begin{center}
        \raisebox{-2.5cm}{\tikzpic{0.5}{
            \draw (0,0) grid (1,5);\draw (1,0) grid (2,1);\draw (2,0) grid (3,5);\draw (3,4) grid (4,5);\draw (4,0) grid (5,5);
            \draw[->, color=gray] (0.5,4.5)--(0.5,0.5);\draw[->, color=gray] (0.5,0.5)--(2.5,0.5);
            \draw[->, color=gray] (2.5,0.5)--(2.5,4.5);\draw[->, color=gray] (2.5,4.5)--(4.5,4.5);
            \draw[->, color=gray] (4.5,4.5)--(4.5,0.5);
        }}
    \end{center}

    Nevertheless, finding the most efficient contraction order is a complex problem~\cite{TN_1}, and we will not delve further into this discussion here.

    \begin{example}\label{ex:cyl}
        In $l \times h$ cylinders.
    \end{example}

    First, consider a \textit{vertically oriented} cylinder, where $l$ is the width and $h$ is the circumference.
    This is equivalent to removing the offset during the assembly of the helicoid or applying a one-step offset in the opposite direction.
    We can achieve this effect using the following Wang tile set:
    \[
        \sum_{\alpha \in \varTheta} \sum_{\beta \in \varTheta} \wtile{$\alpha$}{$\alpha$}{$\beta$}{$\beta$}.
    \]
    Note that we always consider the Wang tile set as a $\bbK$-Wang tile set with weights assigned as either $0$ or $1$, thereby enabling the summation over it.
    By adding a column of positions at the seam (the gray regions in the figure below) where the $l \times h$ chessboard is rolled and joined into a helicoid, assigning this tile to these positions,
    and restricting the characters on the boundary (the bold lines) to be $\sharp$,
    the resulting tilings correspond one-to-one with Wang tilings on the $l \times h$ cylinder.
    The long equal sign in the figure shows the path of these Wang tiles passing horizontal characters.

    \begin{center}
        \raisebox{-1.25cm}{\tikzpic{0.5}{
            \draw (0,0) grid (4,4);
            \node at(2.5,3.5) {$11$};\node at(3.5,3.5) {$21$};\node at(0.5,3.5) {$\ldots$};
            \node at(1.5,3.5) {$l1$};\node at(2.5,2.5) {$12$};\node at(3.5,2.5) {$\ldots$};
            \draw (0,4.5) arc (90:270:0.25);\draw (-0.25,0.25) arc (180:270:0.25);\draw (-0.25,0.25)--(-0.25,4.25);
            \draw (4,0) arc (-90:0:0.25);\draw (4,4) arc (-90:90:0.25);\draw (4.25,0.25)--(4.25,4.25);
            \draw(4,4.5)--(0,4.5);
        }}
        ~~$\Leftrightarrow$
        \raisebox{-1.25cm}{\tikzpic{0.5}{
            \draw (1,0) grid (3,4);\draw (4,1) grid (6,5);
            \node at(4.5,4.5) {$11$};\node at(5.5,4.5) {$21$};\node at(1.5,3.5) {$\ldots$};\node at(2.5,3.5) {$l1$};
            \node at(4.5,3.5) {$12$};\node at(5.5,3.5) {$\ldots$};
            \draw (6.25,4.75) arc (0:90:0.25);\draw (6,0.5) arc (-90:90:0.25);\draw (6.25,4.75)--(6.25,0.75);
            \draw (4,0.5)--(6,0.5);\draw (1,4.5)--(3,4.5);
            \draw (6.25,3.75) arc (0:90:0.25);\draw (6.25,2.75) arc (0:90:0.25);\draw (6.25,1.75) arc (0:90:0.25);
            \draw (1,4.5) arc (90:270:0.25);\draw (0.75,0.25) arc (180:270:0.25);\draw (0.75,0.25)--(0.75,4.25);
            \draw (0.75,1.25) arc (180:270:0.25);\draw (0.75,2.25) arc (180:270:0.25);\draw (0.75,3.25) arc (180:270:0.25);
            \draw[fill=gray!20] (3,0) rectangle (4,5);\draw (3,1)--(4,1);\draw (3,2)--(4,2);\draw (3,3)--(4,3);\draw (3,4)--(4,4);
            \draw (3,0.6)--(4,1.6);\draw (3,0.4)--(4,1.4);\draw (3,1.6)--(4,2.6);\draw (3,1.4)--(4,2.4);
            \draw (3,2.6)--(4,3.6);\draw (3,2.4)--(4,3.4);\draw (3,3.6)--(4,4.6);\draw (3,3.4)--(4,4.4);
            \draw[line width=2pt] (3,4)--(3,5)--(4,5);
            \draw[line width=2pt] (3,0)--(4,0)--(4,1);
        }}
    \end{center}

    Note that the characters on the top and bottom edges of this tile set are horizontal characters, which serve to record and transmit horizontal characters.
    In this case, the size of the state tensor is $\theta^2\sigma^l$, which is only a factor of $\theta$ larger than that of the $l \times h$ chessboard.

    Alternatively, there is another reduction method.
    For convenience, we now consider a \textit{horizontally oriented} cylinder, where $l$ is the width and $h$ is the circumference.
    This is equivalent to gluing the top and bottom edges of the chessboard.

    First, replace each Wang tile in the ($\bbK$-)tile set as follows (keeping the weights unchanged):
    \[
        \wtile{$\alpha$}{$a$}{$b$}{$\beta$} \mapsto \sum_{c \in \varSigma} \wtile{$\alpha$}{$ac$}{$bc$}{$\beta$}.
    \]
    Here, $\varSigma$ is the vertical alphabet of the top edge of the column where the tile set is located
    That is, a new vertical character is appended to the vertical characters on the top and bottom edges, which records the character of the top edge of the column and transmits it downward during the tiling process.

    Thus, when tiling such Wang tiles on the $l \times h$ chessboard, the tiles at the top edge are replaced as follows:
    \[
        \wtile{$\alpha$}{$ac$}{$bc$}{$\beta$} \mapsto \begin{cases} \wtile{$\alpha$}{}{$ba$}{$\beta$}, & a=c; \\
        0, & a \neq c.
        \end{cases}
    \]
    That is, only the cases where the top two characters are equal are retained under the usual constraints (thus, the character $l$ records the top character of this column). 
    Similarly, at the bottom edge, we impose the constraint:
    \[
    \wtile{$\alpha$}{$ac$}{$bc$}{$\beta$} \mapsto \begin{cases} \wtile{$\alpha$}{$ab$}{}{$\beta$}, & b=c; \\
    0, & b \neq c.
    \end{cases}
    \]
    This is equivalent to requiring that the bottom and top characters of the column must both be $l$, effectively gluing the top and bottom edges of the chessboard, as shown in the figure below.

    \begin{center}
        \raisebox{-1.25cm}{\tikzpic{0.5}{
            \draw (0,0) grid (5,5);
            \node at(0.5,4.5) {$11$};\node at(1.5,4.5) {$21$};\node at(3.5,4.5) {$\ldots$};
            \node at(4.5,4.5) {$l1$};\node at(3.5,3.5) {$12$};\node at(4.5,3.5) {$\ldots$};
            \draw (0,5) arc (180:90:0.25);\draw (0,0) arc (180:270:0.25);
            \draw (0.25,5.25)--(5.25,5.25);\draw (0.25,-0.25)--(5.25,-0.25);
            \draw (0,5) arc (180:90:0.25);\draw (0,0) arc (180:270:0.25);\draw (5.5,0)--(5.5,5);
        }}
        ~~$\Leftrightarrow$
        \raisebox{-1.25cm}{\tikzpic{0.5}{
            \draw[fill=gray!20] (0.75,1) rectangle (1.5,4);\draw[fill=gray!20] (2.25,1) rectangle (3,4);
            \draw[fill=gray!20] (3.75,1) rectangle (4.5,4);\draw[fill=gray!20] (5.25,1) rectangle (6,4);
            \draw[fill=gray!20] (6.75,1) rectangle (7.5,4);
            \draw (0,0) rectangle (7.5,5);
            \draw (1.5,0)--(1.5,5);\draw (3,0)--(3,5);\draw (4.5,0)--(4.5,5);\draw (6,0)--(6,5);
            \draw (0,1)--(7.5,1);\draw (0,2)--(7.5,2);\draw (0,3)--(7.5,3);\draw (0,4)--(7.5,4);
            \draw (0.75-0.3,4) arc (180:0:0.3);\draw (0.75+1.5-0.3,4) arc (180:0:0.3);\draw (0.75+3-0.3,4) arc (180:0:0.3);\draw (0.75+4.5-0.3,4) arc (180:0:0.3);\draw (0.75+6-0.3,4) arc (180:0:0.3);
            \draw (0.75-0.5,4) arc (180:0:0.5);\draw (0.75+1.5-0.5,4) arc (180:0:0.5);\draw (0.75+3-0.5,4) arc (180:0:0.5);\draw (0.75+4.5-0.5,4) arc (180:0:0.5);\draw (0.75+6-0.5,4) arc (180:0:0.5);
            \draw (0.75+0.3,4)--(0.75+0.3,1);\draw (0.75+1.5+0.3,4)--(0.75+1.5+0.3,1);\draw (0.75+3+0.3,4)--(0.75+3+0.3,1);\draw (0.75+4.5+0.3,4)--(0.75+4.5+0.3,1);\draw (0.75+6+0.3,4)--(0.75+6+0.3,1);
            \draw (0.75+0.5,4)--(0.75+0.5,1);\draw (0.75+1.5+0.5,4)--(0.75+1.5+0.5,1);\draw (0.75+3+0.5,4)--(0.75+3+0.5,1);\draw (0.75+4.5+0.5,4)--(0.75+4.5+0.5,1);\draw (0.75+6+0.5,4)--(0.75+6+0.5,1);
            \draw (0.75-0.3,1) arc (180:360:0.3);\draw (0.75+1.5-0.3,1) arc (180:360:0.3);\draw (0.75+3-0.3,1) arc (180:360:0.3);\draw (0.75+4.5-0.3,1) arc (180:360:0.3);\draw (0.75+6-0.3,1) arc (180:360:0.3);
            \draw (0.75-0.5,1) arc (180:360:0.5);\draw (0.75+1.5-0.5,1) arc (180:360:0.5);\draw (0.75+3-0.5,1) arc (180:360:0.5);\draw (0.75+4.5-0.5,1) arc (180:360:0.5);\draw (0.75+6-0.5,1) arc (180:360:0.5);
            \node at(0.75,4.5) {$11$};\node at(2.25,4.5) {$21$};\node at(5.25,4.5) {$\ldots$};
            \node at(6.75,4.5) {$l1$};\node at(0.75,3.5) {$12$};\node at(2.25,3.5) {$\ldots$};
        }}
    \end{center}

    Since each vertical character is doubled, the size of the state tensor for the torus case increases from $\theta\sigma^l$ (for the rectangle chessboard case) to $\theta \sigma^{2l}$.
    If such a cylinder is transposed and processed using the first method, the tensor size (note the interchange of horizontal and vertical directions) becomes $\theta^h\sigma^2$.
    Therefore, this approach is only more efficient than the first method when $\sigma^{2l-2} < \theta^{h-1}$.

    \begin{example}
        In $l \times h$ torus.
    \end{example}

    This is equivalent to applying the second method of Example~\ref{ex:cyl} to the vertically oriented cylinder case (processed using the first method of Example~\ref{ex:cyl}) to join the top and bottom edges together.
    In this case, the size of the state tensor becomes $\theta^2\sigma^{2h}$.
    The complexity of the torus case is inevitably much higher than that of all the aforementioned regions with the same $l$ and $h$, which is related to the topological structure of the torus: the torus has no boundary, so regardless of the order of tiling, information about all positions on the edges of the already tiled stack must be recorded.

    Note that the topological structure of torus, like the rectangle chessboard, is transpose-symmetric.
    When $\sigma^{2h-2} > \theta^{2l-2}$, the entire structure can be transposed (the state tensor size becomes $\theta^{2l}\sigma^2$) to enhance computational efficiency.

    \section{Polyomino tilling}

    In this section, we will demonstrate how to transform polyomino tilings in a rectangle chessboard
    (which can be generalized to any area we discussed in the previous chapter) into Wang tilings, thereby enabling the use of our algorithm for enumeration.

    All computational results in this chapter are saved at https://github.com/liang-kai-ruqi,
    and many of them have also been recorded in the corresponding sequences in OEIS.

    Note that we always assume that polyominoes can be rotated, while Wang tiles cannot.

    \subsection{I-shaped polyominoes}

    We first consider the problem of tiling enumeration of I-shaped polyominoes in a rectangle chessboard.
    I-shaped polyominoes refer to rectangular tiles of size $1 \times i$, as shown in the figure below.

    \begin{center}
        \tikzpic{0.6}{
            \drawtile{0}{0};
            \drawdomino{2}{0}{0}{1};
            \node at(2.5,0.5) {$\circlearrowleft_2$};
            \drawtromino{4}{0}{0}{1}{0}{2};
            \node at(4.5,0.5) {$\circlearrowleft_2$};
            \drawtetromino{6}{0}{0}{1}{0}{2}{0}{3};
            \node at(6.5,0.5) {$\circlearrowleft_2$};

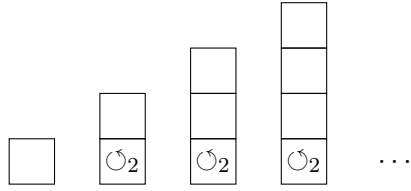
            \node at(8.5,0.5) {$\ldots$};
        }
    \end{center}

    We shall denote the first three by $\Ii, \Iii, \Iiii$ respectively.

    \begin{example}
        Enumerate the tilings in a $l \times h$ chessboard using I-shaped polyominoes.
    \end{example}

    We select the horizontal and vertical alphabets as $\{\sharp, 1\}$, where $\sharp$ is the boundary character.
    The Wang tile set (also treated as a $\bbZ$-tile set with weights $0$ or $1$) is defined as
    \[
        \wtile{$\sharp$}{$\sharp$}{$\sharp$}{$\sharp$}+\wtile{$\sharp$}{$\sharp$}{$\sharp$}{$1$}+\wtile{$1$}{$\sharp$}{$\sharp$}{$\sharp$}
        +\wtile{$\sharp$}{$\sharp$}{$1$}{$\sharp$}+\wtile{$\sharp$}{$1$}{$\sharp$}{$\sharp$}.
    \]
    Here, $\varSigma=\varTheta=\{\sharp, 1\}$, $\sigma = \theta = 2$.
    Hereafter, we shall always assume that the horizontal and vertical alphabets of the Wang tile set consist of all the symbols that appear in the horizontal and vertical directions, respectively.

    Observed that the tiles from the tile set above can and can only form polyominoes with \textit{edges consisting only of boundary characters} as shown below (boundary characters $\sharp$ are no longer drawn).

    \begin{center}
        \tikzpic{0.6}{
            \drawtile{0}{0};
            \drawdomino{2}{0}{0}{1};
            \node at(2.25,1) [right] {$1$};
            \drawtromino{4}{0}{0}{1}{0}{2};
            \node at(4.25,1) [right] {$1$};
            \node at(4.25,2) [right] {$1$};

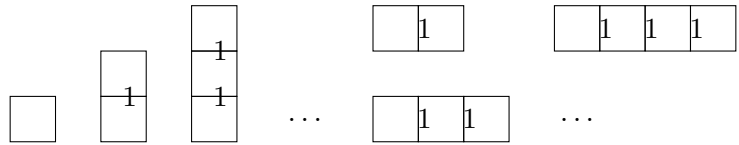
            \node at(6.5,0.5) {$\ldots$};

            \drawdomino{8}{2}{1}{0};
            \node at(8.75,2.5) [right] {$1$};
            \drawtromino{8}{0}{1}{0}{2}{0};
            \node at(8.75,0.5) [right] {$1$};
            \node at(9.75,0.5) [right] {$1$};
            \drawtetromino{12}{2}{1}{0}{2}{0}{3}{0};
            \node at(12.75,2.5) [right] {$1$};
            \node at(13.75,2.5) [right] {$1$};
            \node at(14.75,2.5) [right] {$1$};

            \node at(12.5,0.5) {$\ldots$};
        }
    \end{center}

    This exactly matches the polyominoes required by the problem, and each polyomino corresponds to a \textit{unique} decomposition into the Wang tile set we provided.
    We call such a Wang tile set a \textit{Wang decomposition} of the polyomino.
    Thus, if a polyomino in any region has a Wang decomposition, its tiling enumeration is equal to the tiling enumeration of its Wang decomposition, which our algorithm can be used for.

    For $\bbZ$-Wang tile sets, the above \textit{one-to-one correspondence} can naturally be relaxed to \textit{enumeration correspondence}.
    We will provide a concrete example in Subsection 4 of this section.

    We will provide a more detailed demonstration of the calculation process for this example based on this Wang decomposition,
    as the processes for other examples in this paper are similar.

    In this paper, we always select the data type of integer variables as signed or unsigned 32-bit integer (with a size of 4B).
    Based on the tile set in this example, to ensure no data overflow occurs, the numbers in the state tensor should not exceed $1431655765$.
    Starting from this number, we select $100$ primes downward, obtaining the sequence
    \[
        1431655751, 1431655747, 1431655739, \ldots, 1431653609, 1431653599.
    \]
    The product of these primes is approximately $3.8355 \times 10^{915}$, ensuring it is greater than the number we need to compute.
    This shows that although primes become sparser as numbers grow larger, they are more than sufficient for our algorithm.

    We set the halting mechanism to halt enumerating when the remainder $R$ for each height remains unchanged for $2$ consecutive iterations.
    The estimated probability of misjudgment is less than
    \[
        \frac{1}{1431653609 \times 1431653599} = 2.0496 \times 10^{-18},
    \]
    which is entirely acceptable.
    Moreover, based on the growth trend of the computed arrays, misjudgments can be easily detected.

    We compute all cases for $l \leq 27$ for each $h \leq 27$, and a part of the results are shown in the Table~\ref{tab:I_enum}:
    
    \begin{table}[htbp]
        \centering
        \small
        \begin{tabular}{|r|r|r|r||r|} 
            \hline
            $h$ & $l=1$ & $l=2$ & $l=3$ & $l=h$\\
            \hline
            $1$ &1 &2 &4 &1 \\
            $2$ &2 &7 &29 &7 \\
            $3$ &4 &29 &257 &257 \\
            $4$ &8 &124 &2408 &50128 \\
            $5$ &16 &533 &22873 &50796983 \\
            $6$ &32 &2293 &217969 &264719566561 \\
            $7$ &64 &9866 &2078716 &7063448084710944 \\
            $8$ &128 &42451 &19827701 &963204439792722969647 \\
            $9$ &256 &182657 &189133073 &670733745303300958404439297 \\
            $10$ &512 &785932 &1804125632 &2384351527902618144856749327661056 \\
            $\ldots$ &$\ldots$ &$\ldots$ &$\ldots$ &$\ldots$ \\
            $27$ &\shortstack[r]{\\67\\108\\864}
            &\shortstack[r]{\\46678\\068726\\464909}
            &\shortstack[r]{\\808678\\8352164238\\5015065537}
            &\shortstack[r]{\\112141233\\91782354534780513253590322591761047\\82621009619753629967214626060801782\\86591734152445275752693375509038806\\86290305936932747868007859481561949\\27958742997370483908403590044519371\\84879122372292175736987337417102146\\70845936432144328546887909131495553}\\
            \hline
        \end{tabular}
        \caption{Some results of I-shaped polyominoes tiling enumeration in rectangle chessboards}
        \label{tab:I_enum}
    \end{table}

    It can be seen that, whether for fixed width or when width equals height,
    the number of tiling solutions roughly exhibits exponential growth as the height increases.
    In fact, most polyomino tiling enumeration in a rectangle chessboard follow this pattern,
    while a few exhibit periodicity on top of this, such as periodic fluctuations or periodic occurrences of enumeration being $0$.

    Additionally, the computation time and primes number $n_p$ used for each width are shown in the Table~\ref{tab:I_time}:

    \begin{table}[htbp]
        \centering
        \small

        \begin{tabular}{|c|c|c|c|c|c|c|c|c|}
            \hline
            $l$ & time & $n_p$ & $l$ & time & $n_p$ & $l$ & time & $n_p$ \\
            \hline
            $1$ & 0.002s &3  &$10$ & 0.037s &6 &$19$  & 51.42s &21  \\
            $2$ & 0.003s &4  &$11$ & 0.075s &13 &$20$ & 112.5s &22 \\
            $3$ & 0.002s &5  &$12$ & 0.149s &14 &$21$ & 236.5s &23 \\
            $4$ & 0.003s &6  &$13$ & 0.564s &15 &$22$ & 508.9s &24 \\
            $5$ & 0.003s &7  &$14$ & 0.757s &16 &$23$ & 1101s &25 \\
            $6$ & 0.005s &8  &$15$ & 1.561s &17 &$24$ & 2415s  &26 \\
            $7$ & 0.006s &9  &$16$ & 3.823s &18 &$25$ & 5218s  &27 \\
            $8$ & 0.011s &10 &$17$ & 10.66s &19 &$26$ & 11427s  &28 \\
            $9$ & 0.021s &12 &$18$ & 22.90s &20 &$27$ & 24244s &29 \\
            \hline
        \end{tabular}
        \caption{Calculation time and used primes numbers of I-shaped polyominoes tiling enumeration in rectangle chessboards}
        \label{tab:I_time}
    \end{table}

    The total computation time for the above was 45175s, with $53.667\%$ of the time spent on computing the case with the largest width.
    It can be observed that, for a fixed maximum height, the computation time increases roughly exponentially with the width $l$ of the tiling region:
    Roughly speaking, each increase in width by $1$ doubles the computation time.
    This highlights the significant impact of chessboard width on algorithm complexity.

    Furthermore, the most time-consuming step of the algorithm is always matrix operations (multiplication and modulo),
    while the time spent on other steps (such as solving congruence equations) is negligible.

    Special I-shaped polyominoes can also be decomposed into Wang tiles.
    For example:

    \begin{example}\label{ex:I2+I3}
        Enumerate the tilings in a $l \times h$ chessboard using $\Iii + \Iiii$.
    \end{example}

    These two polyominoes have Wang decompositions (hereafter, characters not shown are all $\sharp$):
    \begin{align*}
        & \wtile{}{}{$2$}{}+\wtile{}{}{}{$2$}+\wtile{}{}{$1$}{}+\wtile{}{}{}{$1$}\\
        +& \wtile{}{$2$}{$1$}{}+\wtile{$2$}{}{}{$1$}+\wtile{}{$1$}{}{}+\wtile{$1$}{}{}{}.
    \end{align*}
    Here, $\sigma = \theta = 3$.

    The tiles in the first row are \textit{root tiles}, while the others, called \textit{branch tiles}, are used for \textit{growing} the tiles.
    The numbers on each tile indicate the number of tiles needed to grow to the right or downward from that position (as shown in the figure below).

    \begin{center}
        \tikzpic{0.6}{
            \drawdomino{0}{0}{0}{1};
            \draw[color=gray,->] (0.5,1.5)--(0.5,0.5);
            \node at(0.25,1) [right] {$1$};
            \drawtromino{2}{0}{0}{1}{0}{2};
            \draw[color=gray,->] (2.5,2.5)--(2.5,1.5);
            \draw[color=gray,->] (2.5,1.5)--(2.5,0.5);
            \node at(2.25,2) [right] {$2$};
            \node at(2.25,1) [right] {$1$};
            \drawdomino{4}{2}{1}{0};
            \draw[color=gray,->] (4.5,2.5)--(5.5,2.5);
            \node at(4.75,2.5) [right] {$1$};
            \drawtromino{4}{0}{1}{0}{2}{0};
            \draw[color=gray,->] (4.5,0.5)--(5.5,0.5);
            \draw[color=gray,->] (5.5,0.5)--(6.5,0.5);
            \node at(4.75,0.5) [right] {$2$};
            \node at(5.75,0.5) [right] {$1$};
        }
    \end{center}

    This approach of \textit{tiles growing} is a general method for designing Wang decompositions of I-shaped polyominoes.
    It is easy to see that the root tiles in the horizontal and vertical directions
    \[
        \wtile{}{}{}{$2$}, \wtile{}{}{$2$}{}
    \]
    determine the total number of tiles,
    while the remaining tiles (branch tiles) are used for growth in various directions.
    If necessary, new \textit{process tiles} of the form
    \[
        \wtile{$(i+1)$}{}{}{$i$}, \wtile{}{$(i+1)$}{$i$}{}
    \]
    can be added.
    Thus, by controlling the possible total number of tiles in each direction using root tiles, this method can construct Wang decompositions for any finite-sized I-shaped polyomino set.
    For example:

    \begin{example}\label{ex:I2}
        Enumerate the tilings in a $l \times h$ chessboard using $\Iii$ (dominoes).
    \end{example}

    It has a Wang decomposition:
    \[
        \wtile{}{}{$1$}{}+\wtile{}{}{}{$1$}+\wtile{}{$1$}{}{}+\wtile{$1$}{}{}{}.
    \]
    Here, $\sigma = \theta = 2$.
    However, this case has a specific formula for the general situation (See Introduction), so we will not compute it further.

    \begin{example}\label{ex:I3}
        Enumerate the tilings in a $l \times h$ chessboard using $\Iiii$ (straight trominoes).
    \end{example}

    It has a Wang decomposition:
    \[
        \wtile{}{}{$2$}{}+\wtile{}{}{}{$2$}
        +\wtile{}{$2$}{$1$}{}+\wtile{$2$}{}{}{$1$}+\wtile{}{$1$}{}{}+\wtile{$1$}{}{}{}.
    \]
    Here, $\sigma = \theta = 3$.

    \subsection{Small polyominoes}

    Clearly, any polyomino tiling (in regions of any shape) can be transformed into Wang tiling by finding a suitable Wang decomposition,
    thereby enabling the use of our algorithm for enumeration.

    In fact,
    \begin{proposition}
        Any finite-area polyomino set must have a Wang decomposition.
    \end{proposition}
    This is because, in the \textit{worst} case, one can enumerate every possible orientation of each polyomino,
    add distinct non-boundary characters to each internal edge, and then separate the tiles.
    It is easy to see that the resulting Wang tile set is a Wang decomposition.
    For example, the Wang decomposition for $\Iii + \Liii$ obtained in this way is
    \begin{center}
        \tikzpic{0.6}{
            \drawdomino{-5}{0}{1}{0};
            \node at(-4.25,0.5) [right] {$1'$};

            \drawdomino{-2}{0}{0}{1};
            \node at(-1.75,1) [right] {$1$};

            \drawtromino{0}{0}{0}{1}{1}{1};
            \node at(0.25,1) [right] {$2$};
            \node at(0.75,1.5) [right] {$2'$};

            \drawtromino{3}{1}{1}{0}{1}{-1};
            \node at(4.25,1) [right] {$3$};
            \node at(3.75,1.5) [right] {$3'$};

            \drawtromino{6}{0}{1}{0}{1}{1};
            \node at(7.25,1) [right] {$4$};
            \node at(6.75,0.5) [right] {$4'$};

            \drawtromino{9}{0}{1}{0}{0}{1};
            \node at(9.25,1) [right] {$5$};
            \node at(9.75,0.5) [right] {$5'$};
        }.
    \end{center}

    However, an excessive number of characters (especially vertical characters) can significantly reduce the efficiency of the algorithm.
    Therefore, for polyominoes with more regular structures,
    it is best to employ some techniques to minimize the number of characters when searching for a Wang decomposition.

    In this subsection, we list the Wang decompositions for all non-empty sets of polyominoes that chosen from the $5$ shapes below.
    \begin{center}
        \tikzpic{0.6}{
            \drawtile{0}{0};

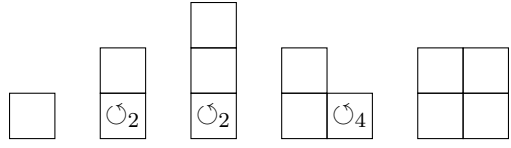
            \drawdomino{2}{0}{0}{1};
            \node at(2.5,0.5) {$\circlearrowleft_2$};

            \drawtromino{4}{0}{0}{1}{0}{2};
            \node at(4.5,0.5) {$\circlearrowleft_2$};

            \drawtromino{6}{0}{0}{1}{1}{0};
            \node at(7.5,0.5) {$\circlearrowleft_4$};

            \drawtetromino{9}{0}{0}{1}{1}{0}{1}{1};
        }
    \end{center}
    We shall denote them by $\Ii, \Iii, \Iiii, \Liii, \Oiv$ respectively.
    Thus, there are $2^5 - 1 = 31$ different non-empty tile sets in total.

    Based on observation, we have:
    \begin{proposition}\label{pro:I1}
        If a ($\bbK$-)Wang tile set $\mathcal{T}$ is a Wang decomposition of a polyomino set $\mathcal{M}$,
        then $\mathcal{T} + \wtile{}{}{}{}$ is a Wang decomposition of $\mathcal{M} + \Ii$.

        (This operation on $\mathcal{T}$, adding \wtile{}{}{}{}, does not change the number of horizontal or vertical characters.)
    \end{proposition}
    Therefore, we will not list polyomino sets that include $\Ii$.

    \begin{example}\label{ex:O4}
        Enumerate the tilings in a $l \times h$ chessboard using $\Oiv$.
    \end{example}

    It has a Wang decomposition:
    \begin{align*}
        \wtile{}{}{$1$}{$2'$}+\wtile{$2'$}{}{$1$}{}+\wtile{}{$1$}{}{}.
    \end{align*}
    Here, $\sigma = \theta = 2$.

    Clearly, these three Wang tiles can only form the shape of $\Oiv$ in the way shown below.
    \begin{center}
        \tikzpic{0.6}{
            \drawtetromino{0}{0}{0}{1}{1}{0}{1}{1};
            \node at(0.75,1.5) [right] {$2'$};
            \node at(0.25,1) [right] {$1$};
            \node at(1.25,1) [right] {$1$};
        }
    \end{center}
    Although there is a vertical line inside the polyomino with a boundary character, it does not matter.

    However, the enumeration of $\Oiv$ tilings on a rectangle chessboard is obvious: it is $1$ if both $l$ and $h$ are even, and $0$ otherwise.
    Thus, we only need to use the above Wang decomposition and Proposition~\ref{pro:I1} to compute the case of $\Ii + \Oiv$.

    \begin{example}
        Enumerate the tilings in a $l \times h$ chessboard using $\Iii + \Oiv$.
    \end{example}

    It has a Wang decomposition:
    \begin{align*}
        \wtile{}{}{$1$}{}+\wtile{}{$1$}{}{}+\wtile{$1$}{}{}{}+\wtile{}{}{}{$1$}
        +\wtile{$2'$}{}{$1$}{}+\wtile{}{}{$1$}{$2'$}.
    \end{align*}
    Here, $\sigma = 2$, $\theta = 3$.

    That is, based on the Wang decomposition of $\Iii$ in Example~\ref{ex:I2}, we add
    \[
        \wtile{$2'$}{}{$1$}{},\wtile{}{}{$1$}{$2'$}.
    \]
    These, together with \wtile{}{$1$}{}{}, exactly form $\Oiv$ (as in the Example~\ref{ex:O4}).

    By analogy with this example, we've observed that:
    \begin{proposition}\label{pro:O4}
        If a ($\bbK$-)Wang tile set $\mathcal{T}$ is a Wang decomposition of a polyomino set $\mathcal{M}$, and it satisfies the following conditions:

        (1) It includes \wtile{}{$1$}{}{} (with weight $1$).

        (2) It does not include the horizontal character $2'$.

        Then,
        \[
            \mathcal{T}+\wtile{}{}{$1$}{$2'$}+\wtile{$2'$}{}{$1$}{}
        \]
        is a Wang decomposition of $\mathcal{M} + \Oiv$.

        (This operation on $\mathcal{T}$ increases the number of horizontal characters by one and does not change the number of vertical characters, having little impact on the algorithm's efficiency.)
    \end{proposition}

    Since all tile sets listed below in this section satisfy the above conditions, we will not list polyomino sets that include $\Oiv$.

    \begin{example}\label{ex:I2+L3}
        Enumerate the tilings in a $l \times h$ chessboard using $\Iii + \Liii$.
    \end{example}

    It has a Wang decomposition:
    \begin{align*}
         & \wtile{}{}{}{$2$}+\wtile{}{}{$2$}{}+\wtile{}{}{$1$}{$1\rtri$}+\wtile{}{}{$1$}{}+\wtile{}{}{}{$1\rtri$}\\
        +& \wtile{$2$}{}{$1$}{}+\wtile{}{$2$}{}{$1\rtri$}+\wtile{$\ltri1$}{$2$}{}{}
        + \wtile{$1\rtri$}{}{}{}+\wtile{}{$1$}{}{}+\wtile{}{}{}{$\ltri1$}.
    \end{align*}
    Here, $\sigma = 3$, $\theta = 4$.

    This is equivalent to selecting the tiles in the first row as \textit{root tiles} for \textit{tree-like growth} (which allows branching), with the following rules:

    (1) The growth direction of the root tiles can only be downward or to the right; the growth direction of other tiles can only be downward, to the left, or to the right.

    (2) A total of $2$ or $3$ tiles are grown (i.e., the sum of the numbers on the root tiles is $1$ or $2$).

    (3) The second growth direction must differ from the first.
    In other words, the corresponding branch tiles (all in the second row) must be \textit{turning}.

    It is easy to verify that this can uniquely obtain each orientation of $\Iii + \Liii$ under rotation (as shown below).
    \begin{center}
        \tikzpic{0.6}{
            \drawdomino{-5}{0}{1}{0};
            \draw[color=gray,->] (-4.5,0.5)--(-3.5,0.5);
            \node at(-4.25,0.5) [right] {$1\rtri$};

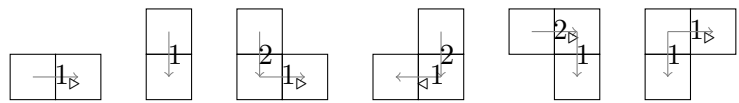
            \drawdomino{-2}{0}{0}{1};
            \draw[color=gray,->] (-1.5,1.5)--(-1.5,0.5);
            \node at(-1.75,1) [right] {$1$};

            \drawtromino{0}{0}{0}{1}{1}{0};
            \draw[color=gray,->] (0.5,1.5)--(0.5,0.5);\node at(0.25,1) [right] {$2$};
            \draw[color=gray,->] (0.5,0.5)--(1.5,0.5);\node at(0.75,0.5) [right] {$1\rtri$};

            \drawtromino{3}{0}{1}{0}{1}{1};
            \draw[color=gray,->] (4.5,1.5)--(4.5,0.5);\node at(4.25,1) [right] {$2$};
            \draw[color=gray,->] (4.5,0.5)--(3.5,0.5);\node at(3.75,0.5) [right] {$\ltri1$};

            \drawtromino{7}{0}{0}{1}{-1}{1};
            \draw[color=gray,->] (6.5,1.5)--(7.5,1.5);\node at(6.75,1.5) [right] {$2\rtri$};
            \draw[color=gray,->] (7.5,1.5)--(7.5,0.5);\node at(7.25,1) [right] {$1$};

            \drawtromino{9}{0}{0}{1}{1}{1};
            \draw[color=gray,->] (9.5,1.5)--(9.5,0.5);\node at(9.75,1.5) [right] {$1\rtri$};
            \draw[color=gray,->] (9.5,1.5)--(10.5,1.5);\node at(9.25,1) [right] {$1$};
        }
    \end{center}
    Note that horizontal growth has two distinct directions, so $\ltri1$ and $1\rtri$ are used to distinguish them.
    Otherwise, an additional polyomino (as shown below) would appear.
    \begin{center}
        \tikzpic{0.6}{
            \drawtetromino{0}{0}{0}{1}{1}{0}{1}{1};
            \node at(0.25,1) [right] {$2$};
            \node at(1.25,1) [right] {$2$};
            \node at(0.75,0.5) [right] {$1$};
        }
    \end{center}
    Thus, by merging the horizontal characters $\ltri1$ and $1\rtri$ into $1$, we obtain the Wang decomposition for $\Iii + \Liii + \Oiv$,
    which has $2$ fewer horizontal characters than the Wang decomposition obtained using Proposition~\ref{pro:O4}.

    \begin{example}
        Enumerate the tilings in a $l \times h$ chessboard using $\Liii$ (right trominoes).
    \end{example}

    By restricting the number of tiles to $3$ in Example~\ref{ex:I2+L3}, i.e., removing the root tiles with $2$ tiles
    \[
        \wtile{}{}{$1$}{}, \wtile{}{}{}{$1\rtri$},
    \]
    the resulting Wang decomposition is
    \begin{align*}
         & \wtile{}{}{}{$2$}+\wtile{}{}{$2$}{}+\wtile{}{}{$1$}{$1\rtri$}\\
        +& \wtile{$2$}{}{$1$}{}+\wtile{}{$2$}{}{$1\rtri$}+\wtile{$\ltri1$}{$2$}{}{}
        +\wtile{$1\rtri$}{}{}{}+\wtile{}{$1$}{}{}+\wtile{}{}{}{$\ltri1$}.
    \end{align*}
    Here, $\sigma = 3$, $\theta = 4$.

    \begin{example}\label{ex:I2+I3+L3}
        Enumerate the tilings in a $l \times h$ chessboard using $\Iii + \Iiii + \Liii$ (dominoes and all trominoes).
    \end{example}

    Simply remove rule (3) from Example~\ref{ex:I2+L3}, i.e., allow tiles to grow without turning.
    The resulting Wang decomposition is
    \begin{align*}
         & \wtile{}{}{}{$2$}+\wtile{}{}{$2$}{}+\wtile{}{}{$1$}{$1\rtri$}+\wtile{}{}{$1$}{}+\wtile{}{}{}{$1\rtri$}\\
        +& \wtile{$2$}{}{$1$}{}+\wtile{$2$}{}{}{$1\rtri$}
        +\wtile{$\ltri1$}{$2$}{}{}+\wtile{}{$2$}{}{$1\rtri$}+\wtile{}{$2$}{$1$}{}
        +\wtile{$1\rtri$}{}{}{}+\wtile{}{$1$}{}{}+\wtile{}{}{}{$\ltri1$}.
    \end{align*}
    Here, $\sigma = 3$, $\theta = 4$.

    \begin{example}\label{ex:I3+L3}
        Enumerate the tilings in a $l \times h$ chessboard using $\Iiii + \Liii$ (all trominoes).
    \end{example}

    Simply restrict the number of tiles in rule (2) of Example~\ref{ex:I2+I3+L3} to $3$, i.e., remove the corresponding root tiles
    \[
        \wtile{}{}{$1$}{}, \wtile{}{}{}{$1\rtri$}.
    \]
    The resulting Wang decomposition is
    \begin{align*}
         & \wtile{}{}{}{$2$}+\wtile{}{}{$2$}{}+\wtile{}{}{$1$}{$1\rtri$}\\
        +& \wtile{$2$}{}{$1$}{}+\wtile{$2$}{}{}{$1\rtri$}
        +\wtile{$\ltri1$}{$2$}{}{}+\wtile{}{$2$}{}{$1\rtri$}+\wtile{}{$2$}{$1$}{}
        +\wtile{$1\rtri$}{}{}{}+\wtile{}{$1$}{}{}+\wtile{}{}{}{$\ltri1$}.
    \end{align*}
    Here, $\sigma = 3$, $\theta = 4$.

    \subsection{Tetrominoes}

    Tetrominoes, also known as $4$-ominoes or \textit{Tetris pieces}, have $7$ different shapes (as shown below).
    \begin{center}
        \tikzpic{0.5}{
            \drawtetromino{0}{0}{1}{0}{0}{1}{0}{2};
            \node at(0.5,0.5) {$\circlearrowleft_4$};
            \node at(1,0) [below] {(L)};
            \drawtetromino{3}{0}{1}{0}{1}{1}{1}{2};
            \node at(3.5,0.5) {$\circlearrowleft_4$};
            \node at(4,0) [below] {(J)};
            \drawtetromino{6}{0}{0}{1}{1}{1}{1}{2};
            \node at(6.5,0.5) {$\circlearrowleft_2$};
            \node at(7,0) [below] {(Z)};
            \drawtetromino{9}{1}{0}{1}{1}{0}{1}{-1};
            \node at(9.5,1.5) {$\circlearrowleft_2$};
            \node at(10,0) [below] {(S)};
            \drawtetromino{12}{0}{0}{1}{0}{2}{1}{1};
            \node at(12.5,0.5) {$\circlearrowleft_4$};
            \node at(13,0) [below] {(T)};
            \drawtetromino{15}{0}{0}{1}{0}{2}{0}{3};
            \node at(15.5,0.5) {$\circlearrowleft_2$};
            \node at(15.5,0) [below] {(I)};
            \drawtetromino{17}{0}{0}{1}{1}{0}{1}{1};
            \node at(18,0) [below] {(O)};
        }
    \end{center}

    \begin{example}
        Enumerate the tilings in a $l \times h$ chessboard using all tetrominoes.
    \end{example}

    Based on the growth rules in Example~\ref{ex:I3+L3} (trominoes), we add one more to the total number of tiles by one, i.e.,

    (1) The growth direction of the root tiles can only be downward or to the right; the growth direction of other tiles can only be downward, to the left, or to the right. Multiple directions can grow simultaneously.

    (2) A total of $4$ tiles are grown (i.e., the sum of the numbers on the root tiles is $3$).

    This results in the following tile set:
    \begin{align*}
         & \wtile{}{}{$3$}{}+\wtile{}{}{$2$}{$1\rtri$}+\wtile{}{}{$1$}{$2\rtri$}\\
        +& \wtile{$\ltri2$}{$3$}{}{}+\wtile{}{$3$}{$2$}{}+\wtile{}{$3$}{}{$2\rtri$}
         +\wtile{$\ltri1$}{$3$}{$1$}{}+\wtile{$\ltri1$}{$3$}{}{$1\rtri$}
         +\wtile{}{$3$}{$1$}{$1\rtri$}\\
        +& \wtile{$\ltri1$}{$2$}{}{}+\wtile{}{$2$}{$1$}{}+\wtile{}{$2$}{}{$1\rtri$}
         +\wtile{$2\rtri$}{}{}{$1\rtri$}+\wtile{$2\rtri$}{}{$1$}{}
         +\wtile{$\ltri1$}{}{}{$\ltri2$}+\wtile{}{}{$1$}{$\ltri2$} \\
        +&\wtile{}{$1$}{}{}+\wtile{$1\rtri$}{}{}{}+\wtile{}{}{}{$\ltri1$}.
    \end{align*}

    The only difference between this Wang tile set and a Wang decomposition of all tetrominoes is that it has 3 different ways of growth into $\Oiv$ (as shown below), causing repeated enumeration of $\Oiv$ three times.
    \begin{center}
        \tikzpic{0.6}{
            \drawtetromino{0}{0}{0}{1}{1}{0}{1}{1};
            \draw[color=gray,->] (0.5,1.5)--(1.5,1.5);
            \draw[color=gray,->] (0.5,1.5)--(0.5,0.5);
            \draw[color=gray,->] (0.5,0.5)--(1.5,0.5);
            \node at(0.25,1) [right] {$2$};
            \node at(0.75,0.5) [right] {$1\rtri$};
            \node at(0.75,1.5) [right] {$1\rtri$};
            \node at(1,0) [below] {($\sqsubset$)};

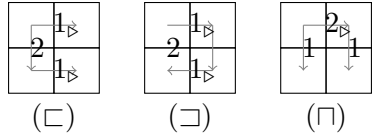
            \drawtetromino{3}{0}{0}{1}{1}{0}{1}{1};
            \draw[color=gray,->] (3.5,1.5)--(4.5,1.5);
            \draw[color=gray,->] (4.5,1.5)--(4.5,0.5);
            \draw[color=gray,->] (4.5,0.5)--(3.5,0.5);
            \node at(3.25,1) [right] {$2$};
            \node at(3.75,0.5) [right] {$1\rtri$};
            \node at(3.75,1.5) [right] {$1\rtri$};
            \node at(4,0) [below] {($\sqsupset$)};

            \drawtetromino{6}{0}{0}{1}{1}{0}{1}{1};
            \draw[color=gray,->] (6.5,1.5)--(6.5,0.5);
            \draw[color=gray,->] (6.5,1.5)--(7.5,1.5);
            \draw[color=gray,->] (7.5,1.5)--(7.5,0.5);
            \node at(6.25,1) [right] {$1$};
            \node at(7.25,1) [right] {$1$};
            \node at(6.75,1.5) [right] {$2\rtri$};
            \node at(7,0) [below] {($\sqcap$)};
        }
    \end{center}
    To resolve this issue, we split the horizontal character $2\rtri$ into $2\rtri$ and $2\rtri'$,
    where $2\rtri'$ indicates that the upstream tile grows both to the right and downward, and $2\rtri$ represents other cases.
    If the tile downstream of $2\rtri'$ grows downward, it forms a $\sqcap$-shape.
    We modify the weight of the downstream tile \wtile{$2\rtri'$}{}{$1$}{} to $-1$,
    so that the weight of $\sqcap$-shaped grown $\Oiv$
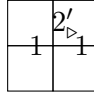
    \begin{center}
        \tikzpic{0.6}{
            \drawtetromino{0}{0}{0}{1}{1}{0}{1}{1};
            \node at(0.75,1.5) [right] {$2\rtri'$};
            \node at(0.25,1) [right] {$1$};
            \node at(1.25,1) [right] {$1$};
        }
    \end{center}
    is $-1$, which exactly cancels out the weights of the other two growth methods to yield $1$, thereby eliminating the repeated enumeration of $\Oiv$.

    The resulting Wang decomposition of tetrominoes (as a $\bbZ$-tile set) is
    \begin{align*}
        & \wtile{}{}{$3$}{}+\wtile{}{}{$2$}{$1\rtri$}+\wtile{}{}{$1$}{$2\rtri'$} \\
        +& \wtile{$\ltri2$}{$3$}{}{}+\wtile{}{$3$}{$2$}{}+\wtile{}{$3$}{}{$2\rtri$}
        +\wtile{$\ltri1$}{$3$}{$1$}{}+\wtile{$\ltri1$}{$3$}{}{$1\rtri$}
        +\wtile{}{$3$}{$1$}{$1\rtri$} \\
        +& \wtile{$\ltri1$}{$2$}{}{}+\wtile{}{$2$}{$1$}{}+\wtile{}{$2$}{$1\rtri$}{}
        +\wtile{$2\rtri$}{}{}{$1\rtri$}+\wtile{$2\rtri$}{}{$1$}{}
        -\wtile{$2\rtri'$}{}{$1$}{}
        +\wtile{$\ltri1$}{}{}{$\ltri2$}+\wtile{}{}{$1$}{$\ltri2$}\\
        +&\wtile{}{$1$}{}{}+\wtile{$1\rtri$}{}{}{}+\wtile{}{}{}{$\ltri1$}.
    \end{align*}
    Here, $\sigma = 4$, $\theta = 6$.

    \subsection{Pentominoes}

    Pentominoes, also known as $5$-ominoes, have $18$ distinct shapes under rotation (as shown below, note that flipped symmetric pentominoes share the same letter as their label):

    \begin{center}
        \tikzpic{0.3}{
            \draw (0,0) rectangle (1,1);\draw (0,1) rectangle (1,2);\draw (0,2) rectangle (1,3);\draw (0,3) rectangle (1,4);\draw (1,0) rectangle (2,1);
            \draw (3,0) rectangle (4,1);\draw (4,0) rectangle (5,1);\draw (4,1) rectangle (5,2);\draw (4,2) rectangle (5,3);\draw (4,3) rectangle (5,4);
            \node at(2.5,0) [below] {(L)};
            \node at(1,1) {$\circlearrowleft_4$};\node at(4,1) {$\circlearrowleft_4$};

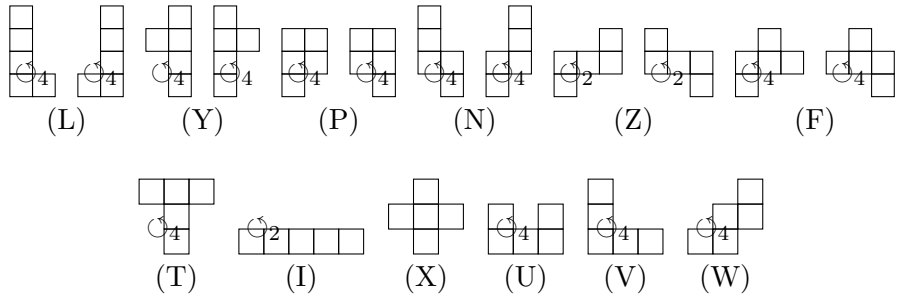
            \draw (6,2) rectangle (7,3);\draw (7,0) rectangle (8,1);\draw (7,1) rectangle (8,2);\draw (7,2) rectangle (8,3);\draw (7,3) rectangle (8,4);
            \draw (11,2) rectangle (10,3);\draw (10,0) rectangle (9,1);\draw (10,1) rectangle (9,2);\draw (10,2) rectangle (9,3);\draw (10,3) rectangle (9,4);
            \node at(8.5,0) [below] {(Y)};
            \node at(7,1) {$\circlearrowleft_4$};\node at(10,1) {$\circlearrowleft_4$};

            \draw (12,0) rectangle (13,1);\draw (12,1) rectangle (13,2);\draw (12,2) rectangle (13,3);\draw (13,1) rectangle (14,2);\draw (13,2) rectangle (14,3);
            \draw (17,0) rectangle (16,1);\draw (17,1) rectangle (16,2);\draw (17,2) rectangle (16,3);\draw (16,1) rectangle (15,2);\draw (16,2) rectangle (15,3);
            \node at(14.5,0) [below] {(P)};
            \node at(13,1) {$\circlearrowleft_4$};\node at(16,1) {$\circlearrowleft_4$};

            \draw (19,1) rectangle (20,2);\draw (18,1) rectangle (19,2);\draw (18,2) rectangle (19,3);\draw (18,3) rectangle (19,4);\draw (19,0) rectangle (20,1);
            \draw (22,1) rectangle (21,2);\draw (23,1) rectangle (22,2);\draw (23,2) rectangle (22,3);\draw (23,3) rectangle (22,4);\draw (22,0) rectangle (21,1);
            \node at(20.5,0) [below] {(N)};
            \node at(19,1) {$\circlearrowleft_4$};\node at(22,1) {$\circlearrowleft_4$};

            \draw (24,0) rectangle (25,1);\draw (24,1) rectangle (25,2);\draw (25,1) rectangle (26,2);\draw (26,1) rectangle (27,2);\draw (26,2) rectangle (27,3);
            \draw (31,0) rectangle (30,1);\draw (31,1) rectangle (30,2);\draw (30,1) rectangle (29,2);\draw (29,1) rectangle (28,2);\draw (29,2) rectangle (28,3);
            \node at(27.5,0) [below] {(Z)};
            \node at(25,1) {$\circlearrowleft_2$};\node at(29,1) {$\circlearrowleft_2$};

            \draw (32,0) rectangle (33,1);\draw (32,1) rectangle (33,2);\draw (33,1) rectangle (34,2);\draw (34,1) rectangle (35,2);\draw (33,2) rectangle (34,3);
            \draw (39,0) rectangle (38,1);\draw (39,1) rectangle (38,2);\draw (38,1) rectangle (37,2);\draw (37,1) rectangle (36,2);\draw (38,2) rectangle (37,3);
            \node at(35.5,0) [below] {(F)};
            \node at(33,1) {$\circlearrowleft_4$};\node at(37,1) {$\circlearrowleft_4$};
        }
    \end{center}
    \begin{center}
        \tikzpic{0.33}{
            \draw (0,2) rectangle (1,3);\draw (1,0) rectangle (2,1);\draw (2,2) rectangle (3,3);\draw (1,1) rectangle (2,2);\draw (1,2) rectangle (2,3);
            \node at(1.5,0) [below] {(T)};
            \node at(1,1) {$\circlearrowleft_4$};

            \draw (4,0) rectangle (5,1);\draw (5,0) rectangle (6,1);\draw (6,0) rectangle (7,1);\draw (7,0) rectangle (8,1);\draw (8,0) rectangle (9,1);
            \node at(6.5,0) [below] {(I)};
            \node at(5,1) {$\circlearrowleft_2$};

            \draw (11,0) rectangle (12,1);\draw (11,1) rectangle (12,2);\draw (10,1) rectangle (11,2);\draw (12,1) rectangle (13,2);\draw (11,2) rectangle (12,3);
            \node at(11.5,0) [below] {(X)};

            \draw (14,0) rectangle (15,1);\draw (14,1) rectangle (15,2);\draw (15,0) rectangle (16,1);\draw (16,0) rectangle (17,1);\draw (16,1) rectangle (17,2);
            \node at(15.5,0) [below] {(U)};
            \node at(15,1) {$\circlearrowleft_4$};

            \draw (18,0) rectangle (19,1);\draw (18,1) rectangle (19,2);\draw (18,2) rectangle (19,3);\draw (19,0) rectangle (20,1);\draw (20,0) rectangle (21,1);
            \node at(19.5,0) [below] {(V)};
            \node at(19,1) {$\circlearrowleft_4$};

            \draw (22,0) rectangle (23,1);\draw (23,0) rectangle (24,1);\draw (23,1) rectangle (24,2);\draw (24,1) rectangle (25,2);\draw (24,2) rectangle (25,3);
            \node at(23.5,0) [below] {(W)};
            \node at(23,1) {$\circlearrowleft_4$};
        }
    \end{center}
    \begin{example}
        Enumerate the tilings in a $l \times h$ chessboard using all pentominoes.
    \end{example}

    We continue using the approach in the previous example, add one more to the total number of tiles.

    On this basis, note that any pentomino containing $\Oiv$ (i.e., the two P-pentominoes, each with $4$ orientations)
    has exactly 3 different ways of growth, such as
    \begin{center}
        \tikzpic{0.6}{
            \drawtromino{0}{0}{0}{1}{0}{2};
            \drawdomino{1}{0}{0}{1};
            \draw[->] (0.5,2.5)--(0.5,1.5);
            \draw[->] (0.5,1.5)--(1.5,1.5);
            \draw[->] (0.5,1.5)--(0.5,0.5);
            \draw[->] (0.5,0.5)--(1.5,0.5);
            \node at(1,0) [below] {($\sqsubset$)};

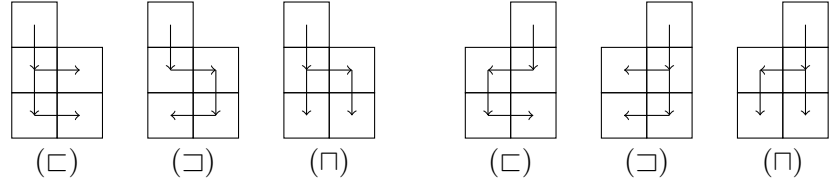
            \drawtromino{3}{0}{0}{1}{0}{2};
            \drawdomino{4}{0}{0}{1};
            \draw[->] (3.5,2.5)--(3.5,1.5);
            \draw[->] (3.5,1.5)--(4.5,1.5);
            \draw[->] (4.5,1.5)--(4.5,0.5);
            \draw[->] (4.5,0.5)--(3.5,0.5);
            \node at(4,0) [below] {($\sqsupset$)};

            \drawtromino{6}{0}{0}{1}{0}{2};
            \drawdomino{7}{0}{0}{1};
            \draw[->] (6.5,2.5)--(6.5,1.5);
            \draw[->] (6.5,1.5)--(7.5,1.5);
            \draw[->] (6.5,1.5)--(6.5,0.5);
            \draw[->] (7.5,1.5)--(7.5,0.5);
            \node at(7,0) [below] {($\sqcap$)};

            \drawtromino{11}{0}{0}{1}{0}{2};
            \drawdomino{10}{0}{0}{1};
            \draw[->] (11.5,2.5)--(11.5,1.5);
            \draw[->] (11.5,1.5)--(10.5,1.5);
            \draw[->] (10.5,1.5)--(10.5,0.5);
            \draw[->] (10.5,0.5)--(11.5,0.5);
            \node at(11,0) [below] {($\sqsubset$)};

            \drawtromino{14}{0}{0}{1}{0}{2};
            \drawdomino{13}{0}{0}{1};
            \draw[->] (14.5,2.5)--(14.5,1.5);
            \draw[->] (14.5,1.5)--(13.5,1.5);
            \draw[->] (14.5,1.5)--(14.5,0.5);
            \draw[->] (14.5,0.5)--(13.5,0.5);
            \node at(14,0) [below] {($\sqsupset$)};

            \drawtromino{17}{0}{0}{1}{0}{2};
            \drawdomino{16}{0}{0}{1};
            \draw[->] (17.5,2.5)--(17.5,1.5);
            \draw[->] (17.5,1.5)--(16.5,1.5);
            \draw[->] (16.5,1.5)--(16.5,0.5);
            \draw[->] (17.5,1.5)--(17.5,0.5);
            \node at(17,0) [below] {($\sqcap$)};
        }.
    \end{center}
    We need to modify the weight of each polyomino containing a $\sqcap$-shaped growth to $-1$.
    This requires splitting the three characters $2\rtri$, $3\rtri$, and $\ltri2$ into two characters each, as in the previous example,
    and then changing the weight of their downstream (either right or left) Wang tiles to $-1$.

    Additionally, since we prohibit upward growth, this tile set cannot form U-shaped pentomino with an upward opening (as shown below).
    \begin{center}
        \tikzpic{0.3}{
            \drawtetromino{0}{0}{0}{1}{1}{0}{2}{0};
            \drawtile{2}{1};
        }
    \end{center}
    Therefore, we add the horizontal character $\mathrm{U}$ and the special tiles
    \[
        \wtile{}{4}{}{$\mathrm{U}$}, \wtile{$\mathrm{U}$}{}{}{$\ltri3$}
    \]
    which are specifically used to form the missing pentomino, as shown below.
    \begin{center}
        \tikzpic{0.6}{
            \drawtetromino{0}{0}{0}{1}{1}{0}{2}{0};
            \drawtile{2}{1};
            \node at(0.25,1) [right] {$4$};
            \node at(2.25,1) [right] {$4$};
            \node at(0.75,0.5) [right] {$\mathrm{U}$};
            \node at(1.75,0.5) [right] {$\ltri3$};
        }
    \end{center}

    This results in the Wang decomposition for all pentominoes:

    \begin{align*}
        & \wtile{}{}{$4$}{}+\wtile{}{}{}{$4\rtri$}+\wtile{}{}{$3$}{$1\rtri$}+\wtile{}{}{$1$}{$3\rtri'$}+\wtile{}{}{$2$}{$2\rtri'$} \\
        + & \wtile{$\ltri3$}{$4$}{}{}+\wtile{}{$4$}{$3$}{}+\wtile{}{$4$}{}{$3\rtri$}+\wtile{$\ltri2'$}{$4$}{$1$}{}+\wtile{$\ltri2$}{$4$}{}{$1\rtri$}+\wtile{$\ltri1$}{$4$}{$2$}{}+\wtile{}{$4$}{$2$}{$1\rtri$}+\wtile{$\ltri1$}{$4$}{$1$}{$1\rtri$}\\
        & ~~ +\wtile{$\ltri1$}{$4$}{}{$2\rtri$}+\wtile{}{$4$}{$1$}{$2\rtri'$}+\wtile{$4\rtri$}{}{$3$}{}+\wtile{$4\rtri$}{}{}{$3\rtri$}+\wtile{$4\rtri$}{}{$2$}{$1\rtri$}+\wtile{$4\rtri$}{}{$1$}{$2\rtri'$}+\wtile{}{$4$}{}{$\mathrm{U}$}\\
        +& \wtile{$\ltri2$}{$3$}{}{}+\wtile{}{$3$}{$2$}{}+\wtile{}{$3$}{}{$2\rtri$}+\wtile{$\ltri1$}{$3$}{$1$}{}+\wtile{$\ltri1$}{$3$}{}{$1\rtri$}+\wtile{}{$3$}{$1$}{$1\rtri$}+\wtile{$3\rtri$}{}{$2$}{}+\wtile{$3\rtri$}{}{}{$2\rtri$}\\
        & ~~ +\wtile{$3\rtri$}{}{$1$}{$1\rtri$}-\wtile{$3\rtri'$}{}{$2$}{}+\wtile{$3\rtri'$}{}{}{$2$}-\wtile{$3'\rtri$}{}{$1$}{$1\rtri$}+\wtile{$\ltri2$}{}{}{$\ltri3$}+\wtile{}{}{$2$}{$\ltri3$}+\wtile{$\ltri1$}{}{$1$}{$\ltri3$}+\wtile{$\mathrm{U}$}{}{}{$\ltri3$}\\
        +&\wtile{$\ltri1$}{$2$}{}{}+\wtile{}{$2$}{$1$}{}+\wtile{}{$2$}{}{$1\rtri$}+\wtile{$2\rtri$}{}{$1$}{}+\wtile{$2\rtri$}{}{}{$1\rtri$}-\wtile{$2\rtri'$}{}{$1$}{}+\wtile{$2\rtri'$}{}{}{$1\rtri$} \\
        & ~~ +\wtile{$\ltri1$}{}{}{$\ltri2$}+\wtile{}{}{$1$}{$\ltri2$}+\wtile{$\ltri1$}{}{}{$\ltri2'$}-\wtile{}{}{$1$}{$\ltri2'$} \\
        +& \wtile{}{$1$}{}{}+\wtile{$\ltri1$}{}{}{}+\wtile{}{}{}{$1\rtri$}.
    \end{align*}
    Here, $\sigma=5$, $\theta=12$.

    Clearly, in this case, the enumeration is positive only when at least one of $l$ or $h$ (without loss of generality, assume it's $l$) is divisible by $5$.
    Table~\ref{tab:5omino_enum} below presents the calculation results for the cases where $l = 5, 10$ (while $h$ ranges over $1,2,\ldots,20$), with computation times of 0.4228s and 3789s.
    Our algorithm achieved the first successful computation for the cases of $l=10$ with large $h$.

    \begin{table}[htbp]
        \centering
        \small

        \begin{tabular}{|r|r|r|}
            \hline
            $h=$ & $l=5$ & $l=10$ \\
            \hline
            $1$  &1  &1  \\
            $2$  &5  &45  \\
            $3$  &56  &7670  \\
            $4$  &501  &890989  \\
            $5$  &4006  &101698212  \\
            $6$  &27950  &7845888732  \\
            $7$  &214689  &756605877809  \\
            $8$  &1696781  &75996685446347  \\
            $9$  &13205354  &7470920047174798  \\
            $10$  &101698212  &729748655181974778  \\
            $11$  &782267786  &70521242596066128006  \\
            $12$  &6048166230  &6882943628424155149082  \\
            $13$  &46799177380  &672858933871350579734838  \\
            $14$  &361683136647  &65670854176387745944044415  \\
            $15$  &2793722300087  &6406383348267533424844337077  \\
            $16$  &21583392631817  &624874119278590450628206097405  \\
            $17$  &166790059833039  &60978146945443555311094030206323  \\
            $18$  &1288885349447958  &5950711244486170431626902119957082  \\
            $19$  &9959188643348952  &580653334431250399171093742069162662  \\
            $20$  &76953117224941654  &56657284915840468039405015713225758536  \\
            \hline
        \end{tabular}
        \caption{Some results of pentominoes tiling enumeration in rectangle chessboards}
        \label{tab:5omino_enum}
    \end{table}

    \subsection{Related OEIS sequences}

    Table~\ref{tab:OEIS_enum} lists all the polyomino tiling enumeration presented in this paper that have corresponding sequences recorded in the OEIS, as well as the precious range and the range of results calculated by our algorithm
    (utilizing transpose symmetry, as well as certain special cases where the enumeration is $0$).
    The device used for calculation is a conventional computer (using a single core of AMD ryzen 7 7840hs CPU).

    \begin{table}[htbp]
        \centering
        \small

        \begin{tabular}{|c|c|c|c|c|}
            \hline
            \shortstack{\\polyomino\\set} & \shortstack{\\ OEIS sequences and\\ precious range} & \shortstack{\\calculation\\range} & time & \shortstack{\\ maximum\\ enumeration} \\
            \hline
            I-shaped & \shortstack{\\ A254414:~$l+h\leq 30$\\ A254127:~$l=h\leq 15$}  & $l+h\leq 55$ & \shortstack{\\total:~53795s\\$l=27$:~28489s} & \shortstack{\\ $3.306\times 10^{262}$\\ $(l=27, h=28)$}\\
            \hline
            $\Iii+\Iiii$ & \shortstack{\\ A219866:~$l+h\leq 22$\\ A219874:~$l=h\leq 11$} & $l+h\leq 35$ & \shortstack{\\total:~9468s\\$l=17$:~6362s} & \shortstack{\\ $ 7.757\times 10^{59}$\\ $(l=17 ,h=18)$} \\
            \hline
            $\Liii$ & A351322:~$l+\frac{h}{3}\leq 31 $ & $l+\frac{h}{3}\leq 35$ & \shortstack{\\total:~62954s\\$l=17$:~42719s} & \shortstack{\\ $3.307 \times 10^{94}$\\ $(l=17, h=54)$} \\
            \hline
            $\Ii+\Liii$ & \shortstack{\\ A220054:~$l+h\leq 29$\\ A220061:~$l=h\leq 16$} &  $l+h\leq 35$ & \shortstack{\\total:~17955s\\$l=17$:~12106s} &\shortstack{\\ $ 3.046\times 10^{77}$\\ $(l=17 ,h=18)$} \\
            \hline
            $\Iii+\Liii$ & \shortstack{\\ A219987:~$l+h\leq 21$\\ A219994:~$l=h\leq 17$} &  $l+h\leq 35$ & \shortstack{\\total:~16380s\\$l=17$:~10987s} &\shortstack{\\ $ 1.167\times 10^{68}$\\ $(l=17 ,h=18)$} \\
            \hline
            $\Ii+\Iii+\Liii$ & \shortstack{\\ A353877:~$l, h\leq 9$\\ A353934:~$l=h\leq 16$} & $l+h\leq 35$ & \shortstack{\\total:~22518s\\$l=17$:~15107s} &\shortstack{\\ $1.813 \times 10^{107}$\\ $(l=17 ,h=18)$} \\
            \hline
            $\Iiii+\Liii$ &  A233320:~$l+h\leq 28$ &  $l+h\leq 35$ & \shortstack{\\total:~13299s\\$l=17$:~8928s} &\shortstack{\\ $2.173 \times 10^{49}$\\ $(l=17,h=18)$} \\
            \hline
            $\Ii+\Iiii+\Liii$ & \shortstack{\\ A270061:~$l+h\leq 25$\\ A270071:~$l=h\leq 12$} & $l+h\leq 35$& \shortstack{\\total:~19344s\\$l=17$:~13017s} &\shortstack{\\ $6.584\times 10^{89}$\\ $(l=17,h=18)$} \\
            \hline
            $\Iii+\Iiii+\Liii$ & \shortstack{\\ A364457:~$l+h\leq 25$\\ A364504:~$l=h\leq 12$} & $l+h\leq 35$ & \shortstack{\\total:~18331s\\$l=17$:~12279s} &\shortstack{\\ $ 4.358\times 10^{79}$\\ $(l=17,h=18)$} \\
            \hline
            $\Ii+\Oiv$ & \shortstack{\\ A245013:~$l+h\leq 45$\\ A063443:~$l=h\leq 40$} & $l+h\leq 55$ & \shortstack{\\total:~19667s\\$l=27$:~9769s} &\shortstack{\\ $8.952 \times 10^{92}$\\ $(l=27,h=28)$} \\
            \hline
            $\Iii+\Oiv$ & A352431:~$l+\frac{h}{2}\leq 40$ & $l+\frac{h}{2}\leq 53$ & \shortstack{\\total:~81727s\\$l=26$:~41408s} &\shortstack{\\ $1.644 \times 10^{231}$\\ $(l=26,h=54)$} \\
            \hline
            $\Ii+\Iii+\Oiv$ & \shortstack{\\ A352589:~$l,h\leq 9$\\ A353777,$l=h\leq 17$} &  $l+h\leq 55$ & \shortstack{\\total:~75583s\\$l=27$:~38039s} &\shortstack{\\ $1.217 \times 10^{224}$\\ $(l=27,h=28)$} \\
            \hline
            $\Iiii+\Oiv$ & \shortstack{\\ A219967:~$l+h\leq 27$\\ A219975:~$l=h\leq 12$} &  $l+h\leq 35$ & \shortstack{\\total:~8339s\\$l=17$:~5589s} &\shortstack{\\ $9.109 \times 10^{27}$\\ $(l=17,h=18)$} \\
            \hline
            $\Liii+\Oiv$ & \shortstack{\\ A219946:~$l+h\leq 35$\\ A219952:~$l=h\leq 17$} &  $l+h\leq 35$ & \shortstack{\\total:~14942s\\$l=17$:~10017s} &\shortstack{\\ $1.467 \times 10^{36}$\\ $(l=17,h=18)$} \\
            \hline
            $\Ii+\Liii+\Oiv$ & \shortstack{\\ A353963:~$l,h\leq 9$\\ A354067:~$l=h\leq 16$} &  $l+h\leq 35$ & \shortstack{\\total:~26400s\\$l=17$:~17712s} &\shortstack{\\ $1.103 \times 10^{81}$\\ $(l=17,h=18)$} \\
            \hline
            tetrominoes & \shortstack{\\ A230031:~$l+h\leq 20$\\ A263425:~$l=h\leq 10$} &  $l+h\leq 27$ & \shortstack{\\total:~8879s\\$l=12$:~6494s} &\shortstack{\\ $8.980 \times 10^{31}$\\ $(l=12,h=15)$} \\
            \hline
            pentominoes &  A233427:~$l+h\leq 17$ &  $l+h\leq 26$ & \shortstack{\\total:~16814s\\$l=11$:~13595s} &\shortstack{\\$7.977 \times 10^{30}$\\ $(l=11,h=15 )$} \\
            \hline
            \shortstack{\\ $\Ii~+$\\pentominoes} &  A278657:~$l+h\leq 15$ &  $l+h\leq 23$ & \shortstack{\\total:~17578s\\$l=11$:~14482s} &\shortstack{\\ $ 4.576\times 10^{38}$\\ $(l=11 ,h=12 )$} \\
            \hline
        \end{tabular}
        \caption{OEIS sequences related to polyomino tiling enumeration and calculation of our algorithm}
        \label{tab:OEIS_enum}
    \end{table}

    Below the total calculation time is for the case of calculating the maximum width,
    and the ratio of it to the total time consumption can roughly reflect the rate at which the total time consumption increases with the width.
    The rightmost column shows the approximate value of the largest count in the calculation results (although we calculate the exact value).

    It can be observed that, compared to the traditional algorithms used in the original records,
    our algorithm improves computational efficiency in the vast majority of cases, thereby significantly expanding the computable range,
    especially for the cases with larger enumeration and more complex polyomino sets.

    \section{Further generalizations and open problems}

    \subsection{Variants of tile shape}

    Our algorithm can also be generalized to the enumeration of tilings of more types.
    For example, \textit{polyiamonds} and \textit{polyhexes} are variants of polyominoes, where the unit tiles are equilateral triangles and regular hexagons, respectively, instead of squares.
    All the shapes of $4$-iamonds ($3$ different shapes) and $3$-hexes ($3$ different shapes) are shown in the figures below.

    \begin{center}
        \tikzpic{0.6}{
            \draw (0,0)--(2,0)--(1,1.732)--(0,0);\draw (1,0)--(1.5,0.866)--(0.5,0.866)--(1,0);
            \draw (3,0)--(5,0)--(5.5,0.866)--(3.5,0.866)--(4,0)--(4.5,0.866)--(5,0);\draw (3,0)--(3.5,0.866);
            \draw (6.5,0.866)--(8.5,0.866)--(8,0)--(7,0)--(6.5,0.866)--(7,1.732)--(8,0);\draw (7.5,0.866)--(7,0);
        }
        ~~~~~~~~
        \tikzpic{0.3}{
            \draw (0,0)--(1*0.866,-0.5)--(2*0.866,0)--(2*0.866,1)--(1*0.866,1.5)--(0,1)--(0,0);
            \draw (2*0.866,0)--(3*0.866,-0.5)--(4*0.866,0)--(4*0.866,1)--(3*0.866,1.5)--(2*0.866,1)--(2*0.866,0);
            \draw (4*0.866,0)--(5*0.866,-0.5)--(6*0.866,0)--(6*0.866,1)--(5*0.866,1.5)--(4*0.866,1)--(4*0.866,0);

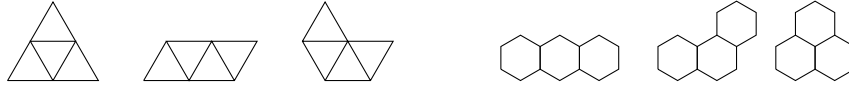
            \draw (8*0.866,0)--(9*0.866,-0.5)--(10*0.866,0)--(10*0.866,1)--(9*0.866,1.5)--(8*0.866,1)--(8*0.866,0);
            \draw (10*0.866,0)--(11*0.866,-0.5)--(12*0.866,0)--(12*0.866,1)--(11*0.866,1.5)--(10*0.866,1)--(10*0.866,0);
            \draw (12*0.866,1)--(13*0.866,1.5)--(13*0.866,2.5)--(12*0.866,3)--(11*0.866,2.5)--(11*0.866,1.5);

            \draw (14*0.866,0)--(15*0.866,-0.5)--(16*0.866,0)--(16*0.866,1)--(15*0.866,1.5)--(14*0.866,1)--(14*0.866,0);
            \draw (16*0.866,0)--(17*0.866,-0.5)--(18*0.866,0)--(18*0.866,1)--(17*0.866,1.5)--(16*0.866,1)--(16*0.866,0);
            \draw (17*0.866,1.5)--(17*0.866,2.5)--(16*0.866,3)--(15*0.866,2.5)--(15*0.866,1.5);
        }
    \end{center}

    The enumeration problems for these variants can also be transformed into similar tensor operations using the techniques of \textit{helicoids} and \textit{Wang partitions} introduced in this paper, and then computed using modular enumeration for their (weighted) enumeration.
    The only difference is that these variants introduce an additional directional alphabet, and the rules for index rotation in matrix operations need to be redesigned to ensure that the first few indices always participate in the multiplication.
    The details are straightforward and can be inferred from the movement of the initial point defining the index order, as shown in the figures below.

    \begin{center}
        \raisebox{-1.5cm}{\tikzpic{0.6}{
            \draw (0,1)--(0.5,1)--(1,2)--(5,2)--(5.25,2.5)--(0.75,2.5)--(0,1);
            \draw (2.5,2)--(2.5,1.5);\draw (3.5,2)--(3.5,1.5);\draw (4.5,2)--(4.5,1.5);
            \node at(2.5,1.5) [below] {$a_1$};\node at(3.5,1.5) [below] {$\ldots$};\node at(4.5,1.5) [below] {$a_l$};
            \draw (1,0)--(2,0)--(1.5,1)--(1,0);
            \draw (1.25,0.5)--(0.75,1.5);\node at (1,1) {$\alpha_1$};
            \draw (1.75,0.5)--(2.25,0.5);\node at(2.25,0.5) [right] {$\gamma_1$};
            \draw (1.5,0)--(1.5,-0.5);\node at(1.5,-0.5) [below] {$b_1$};
            \node at(0.5,1) {$\circ$};\draw[->] (0.5,1)--(0.65,1.3);
        }}
        $\rightarrow$
        \raisebox{-1cm}{\tikzpic{0.6}{
            \draw (0,1)--(1.5,1)--(1,2)--(5,2)--(5.25,2.5)--(0.75,2.5)--(0,1);
            \draw (3.5,2)--(3.5,1.5);\draw (4.5,2)--(4.5,1.5);
            \node at(3.5,1.5) {$\ldots$};\node at(4.5,1.5) [below] {$a_l$};
            \draw (2.75,0.5)--(1.25,1.5);\draw (3,1)--(1.5,2);
            \node at(2.5,1.5) {$a_1$};\node at(2,1) {$\gamma_1$};
            \draw (3,0)--(3.5,1)--(2.5,1)--(3,0);
            \node at(1.5,1) {$\circ$};\draw[->] (1.5,1)--(1.35,1.3);
            \draw[->] (0.5,1.2) arc (110:70:1);
            \draw (3.25,0.5)--(3.75,0.5);\node at(3.75,0.5) [right] {$\alpha_2$};
            \draw (0.5,1)--(0.5,0.5);\node at(0.5,0.5) [below] {$b_1$};
        }}
        $\rightarrow$
        \raisebox{-0.5cm}{\tikzpic{0.6}{
            \draw (0,1)--(0.5,1)--(1,2)--(6,2)--(6.25,2.5)--(0.75,2.5)--(0,1);
            \draw (2.5,2)--(2.5,1.5);\draw (3.5,2)--(3.5,1.5);\draw (4.5,2)--(4.5,1.5);\draw (5.5,2)--(5.5,1.5);
            \node at(2.5,1.5) [below] {$a_2$};\node at(3.5,1.5) [below] {$\ldots$};\node at(4.5,1.5) [below] {$a_l$};\node at(5.5,1.5) [below] {$b_1$};
            \node at(0.5,1) {$\circ$};\draw[->] (0.5,1)--(0.65,1.3);
            \draw (0.75,1.5)--(1.25,1.5);\node at(1.25,1.5) [right] {$\alpha_2$};
        }}
    \end{center}

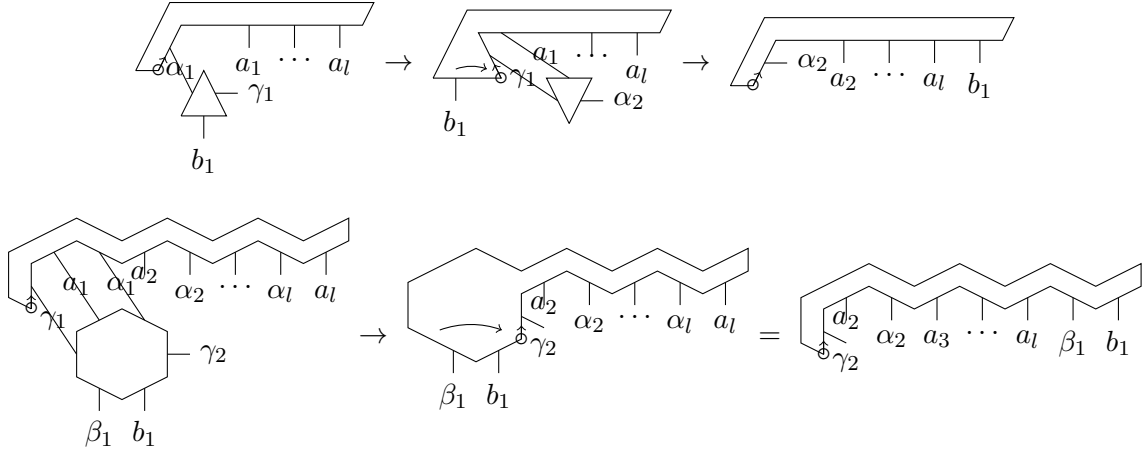
    \begin{center}
        \raisebox{-1.5cm}{\tikzpic{0.6}{
            \draw (1,1)--(1,2)--(2,2.5)--(3,2)--(4,2.5)--(5,2)--(6,2.5)--(7,2)--(8,2.5)--(8,3)--(7,2.5)--(6,3)--(5,2.5)--(4,3)--(3,2.5)--(2,3)--(0.5,2.25)--(0.5,1.25)--(1,1);
            \draw (3.5,2.25)--(3.5,1.75); \node at(3.5,1.75) {$a_2$};
            \draw (4.5,2.25)--(4.5,1.75); \node at(4.5,1.75) [below] {$\alpha_2$};
            \draw (5.5,2.25)--(5.5,1.75); \node at(5.5,1.75) [below] {$\ldots$};
            \draw (6.5,2.25)--(6.5,1.75); \node at(6.5,1.75) [below] {$\alpha_l$};
            \draw (7.5,2.25)--(7.5,1.75); \node at(7.5,1.75) [below] {$a_l$};
            \draw (2,0.5)--(3,1)--(4,0.5)--(4,-0.5)--(3,-1)--(2,-0.5)--(2,0.5);
            \draw (2,0)--(1,1.5);\node at(1.5,0.75) {$\gamma_1$};
            \draw (2.5,0.75)--(1.5,2.25);\node at(2,1.5) {$a_1$};
            \draw (3.5,0.75)--(2.5,2.25);\node at(3,1.5) {$\alpha_1$};
            \draw (4,0)--(4.5,0); \node at(4.5,0) [right] {$\gamma_2$};
            \draw (2.5,-0.75)--(2.5,-1.25); \node at (2.5,-1.25) [below] {$\beta_1$};
            \draw (3.5,-0.75)--(3.5,-1.25); \node at (3.5,-1.25) [below] {$b_1$};
            \node at(1,1) {$\circ$}; \draw[->] (1,1)--(1,1.3);
        }}
        $\rightarrow$
        \raisebox{-1cm}{\tikzpic{0.6}{
            \draw (1,1)--(2,0.5)--(3,1)--(3,2)--(4,2.5)--(5,2)--(6,2.5)--(7,2)--(8,2.5)--(8,3)--(7,2.5)--(6,3)--(5,2.5)--(4,3)--(3,2.5)--(2,3)--(0.5,2.25)--(0.5,1.25)--(1,1);
            \draw (3,1.5)--(3.5,1.25); \node at(3.5,1.25) [below] {$\gamma_2$};
            \draw (3.5,2.25)--(3.5,1.75); \node at(3.5,1.75) {$a_2$};
            \draw (4.5,2.25)--(4.5,1.75); \node at(4.5,1.75) [below] {$\alpha_2$};
            \draw (5.5,2.25)--(5.5,1.75); \node at(5.5,1.75) [below] {$\ldots$};
            \draw (6.5,2.25)--(6.5,1.75); \node at(6.5,1.75) [below] {$\alpha_l$};
            \draw (7.5,2.25)--(7.5,1.75); \node at(7.5,1.75) [below] {$a_l$};
            \draw (1.5,0.75)--(1.5,0.25); \node at (1.5,0.25) [below] {$\beta_1$};
            \draw (2.5,0.75)--(2.5,0.25); \node at (2.5,0.25) [below] {$b_1$};
            \node at(3,1) {$\circ$}; \draw[->] (3,1)--(3,1.3);
            \draw[->] (1.2,1.2) arc (110:70:2);
        }}
        $=$
        \raisebox{-0.5cm}{\tikzpic{0.6}{
            \draw (1,1)--(1,2)--(2,2.5)--(3,2)--(4,2.5)--(5,2)--(6,2.5)--(7,2)--(8,2.5)--(8,3)--(7,2.5)--(6,3)--(5,2.5)--(4,3)--(3,2.5)--(2,3)--(0.5,2.25)--(0.5,1.25)--(1,1);
            \draw (1.5,2.25)--(1.5,1.75); \node at(1.5,1.75) {$a_2$};
            \draw (2.5,2.25)--(2.5,1.75); \node at(2.5,1.75) [below] {$\alpha_2$};
            \draw (3.5,2.25)--(3.5,1.75); \node at(3.5,1.75) [below] {$a_3$};
            \draw (4.5,2.25)--(4.5,1.75); \node at(4.5,1.75) [below] {$\ldots$};
            \draw (5.5,2.25)--(5.5,1.75); \node at(5.5,1.75) [below] {$a_l$};
            \draw (6.5,2.25)--(6.5,1.75); \node at(6.5,1.75) [below] {$\beta_1$};
            \draw (7.5,2.25)--(7.5,1.75); \node at(7.5,1.75) [below] {$b_1$};
            \draw (1,1.5)--(1.5,1.25); \node at(1.5,1.25) [below] {$\gamma_2$};

            \node at(1,1) {$\circ$}; \draw[->] (1,1)--(1,1.3);

        }}
    \end{center}

    Similarly, our algorithm can be generalized and applied to any variant of Wang tiles that exhibits periodicity, as well as the tile sets that can be decomposed into such Wang tiles.

    \subsection{Tilings with restrictions}

    In certain tiling enumeration problems, constraints may be imposed on the number of occurrences of specific tiles or require the tilings to possess certain symmetries.
    We provide two sufficiently general examples corresponding to these two types of constraints, respectively, to demonstrate how our algorithm can be applied to handle these restricted cases.

    \begin{example}
        Enumerate the ways to place $\Iii$ and $\Iiii$ in a $l \times h$ chessboard, such that the total number of tiles is n.
    \end{example}

    Clearly, the placement of $\Iii+\Iiii$ in a chessboard corresponds one-to-one with the tiling of $\Ii+\Iii+\Iiii$ on the same chessboard, as it suffices to fill all the remaining empty spaces using $\Ii$.
    As for the restriction on the number of tiles, it is noteworthy that during tiling in a helicoid, the horizontal characters precisely scan all positions.
    Therefore, the horizontal characters can be utilized to record information related to the number of tiles.

    We use $\{0, 1, \ldots, \nu-1\}\times\varTheta$ to replace the original horizontal alphabet $\varTheta$, where $\nu$ is greater than all the values of $n$ that we need to calculate.
    In this specific case, since each Wang partition of $\Iii$ or $\Iiii$ has exactly one root tile, the number of tiles on the chessboard is precisely equal to the number of root tiles.
    Therefore, we perform the following substitution:
    \[
        \wtile{$\alpha$}{$a$}{$b$}{$\beta$}\mapsto
        \begin{cases}
            \displaystyle\sum_{i=0}^{\nu-2}\wtile{$i\alpha$}{$a$}{$b$}{$(i+1)\beta$}, & \mathrm{root~tile};\\
            \displaystyle\sum_{i=0}^{\nu-1}\wtile{$i\alpha$}{$a$}{$b$}{$i\beta$}, & \mathrm{not~root~tile}.
        \end{cases}
    \]
    Then, by replacing $\sharp\sharp\ldots\sharp$ with $0\sharp\sharp\ldots\sharp$ as the initial word for tiling, the resulting tensor will encompass the enumeration of all tilings that restrict the number of tiles with $n = 0, 1, \ldots, \nu-1$,
    as the weight of $0\sharp\sharp\ldots\sharp$, $1\sharp\sharp\ldots\sharp$, $\ldots$, $(\nu-1)\sharp\sharp\ldots\sharp$, precisely corresponding to the leading $\nu$ terms of the tensor.
    Compared to the case without the number restriction, the space complexity increased by a factor of $\nu$ in this case.

    Furthermore, if the number restriction is modulo $\nu$ with remainder $n$, simply modify the substitution described above as follows:
    \[
        \wtile{$\alpha$}{$a$}{$b$}{$\beta$}\mapsto
        \begin{cases}
            \wtile{$(\mu-1)\alpha$}{$a$}{$b$}{$0\beta$}+\displaystyle\sum_{i=0}^{\nu-2}\wtile{$i\alpha$}{$a$}{$b$}{$(i+1)\beta$}, & \mathrm{root~tile};\\
            \displaystyle\sum_{i=0}^{\nu-1}\wtile{$i\alpha$}{$a$}{$b$}{$i\beta$}, & \mathrm{not~root~tile}.
        \end{cases}
    \]

    \begin{example}
        Enumerate the ways to place $\Iii+\Iiii$ in a $l \times h$ chessboard that the placements are left-right symmetric.
    \end{example}

    We still regard it as a tiling of $\Ii +\Iii+\Iiii$.
    To preserve symmetry, we fold the chessboard along its vertical axis, as shown below.
    \begin{center}
        \raisebox{-1cm}{\tikzpic{0.5}{
            \draw (0,0) grid (6,4);
        }}
        $\rightarrow$
        \raisebox{-1cm}{\tikzpic{0.5}{
            \draw (0,0) grid (3,4);
            \draw (0.2,-0.2)--(0.2,3.8);\draw (1.2,-0.2)--(1.2,3.8);\draw (2.2,-0.2)--(2.2,3.8);
            \draw (0.2,3.8)--(1.2,3.8)--(2.2,3.8)--(3,4);\draw (0.2,2.8)--(1.2,2.8)--(2.2,2.8)--(3,3);
            \draw (0.2,1.8)--(1.2,1.8)--(2.2,1.8)--(3,2);\draw (0.2,0.8)--(1.2,0.8)--(2.2,0.8)--(3,1);
            \draw (0.2,-0.2)--(1.2,-0.2)--(2.2,-0.2)--(3,0);
        }}~~~~
        \raisebox{-1cm}{\tikzpic{0.5}{
            \draw (0,0) grid (5,4);
        }}
        $\rightarrow$
        \raisebox{-1cm}{\tikzpic{0.5}{
            \draw (0,0) grid (2.5,4);
            \draw (0.2,-0.2)--(0.2,3.8);\draw (1.2,-0.2)--(1.2,3.8);\draw (2.2,-0.2)--(2.2,3.8);
            \draw (0.2,3.8)--(1.2,3.8)--(2.2,3.8)--(2.5,4);\draw (0.2,2.8)--(1.2,2.8)--(2.2,2.8)--(2.5,3);
            \draw (0.2,1.8)--(1.2,1.8)--(2.2,1.8)--(2.5,2);\draw (0.2,0.8)--(1.2,0.8)--(2.2,0.8)--(2.5,1);
            \draw (0.2,-0.2)--(1.2,-0.2)--(2.2,-0.2)--(2.5,0);
            \draw[dashed] (2.5,0)--(2.5,4);
        }}
    \end{center}

    When the width of the chessboard is even, each position after folding consists of two stacked tiles.
    We restrict these two tiles to originate from corresponding positions in two left-right symmetric polyominoes.
    That is, We obtain the ($\bbK$-) tile set by
    \[
        \wtile{$\alpha$}{$a$}{$b$}{$\beta$}, \wtile{$\alpha'$}{$a'$}{$b'$}{$\beta'$} \mapsto \wtile{$\alpha\beta'$}{$aa'$}{$bb'$}{$\beta\alpha'$},
    \]
    where \wtile{$\alpha$}{$a$}{$b$}{$\beta$} and \wtile{$\alpha'$}{$a'$}{$b'$}{$\beta'$} range over all Wang tiles that may appear in the left-right symmetric positions within a single polyomino or across two polyominoes.
    For the Wang decomposition of $\Ii +\Iii+\Iiii$ (see Example~\ref{ex:I2+I3}), it is
    \begin{gather*}
        \wtile{}{}{}{$2$}\leftrightarrow\wtile{$1$}{}{}{},~~\wtile{}{}{}{$1$}\leftrightarrow\wtile{$1$}{}{}{},\\
        \wtile{$2$}{}{}{$1$}\!\!\hookleftarrow,~~\wtile{}{}{$2$}{}\!\!\hookleftarrow,~~ \wtile{}{$2$}{$1$}{}\!\!\hookleftarrow,~~
        \wtile{}{$1$}{}{}\!\!\hookleftarrow,~~\wtile{}{}{$1$}{}\!\!\hookleftarrow,~~\wtile{}{}{}{}\!\!\hookleftarrow.
    \end{gather*}
    The tiles in the second row are self-corresponding, meaning they only appear in positions that are symmetric with respect to the left-right symmetry within a single polyomino.
    Thus, the resulting tile set is obtained as follows:
    \begin{align*}
        & \wtile{}{}{}{$21$}+\wtile{$12$}{}{}{}+\wtile{}{}{}{$11$}+\wtile{$11$}{}{}{}\\
        + & \wtile{$21$}{}{}{$12$}+\wtile{}{}{$22$}{}+\wtile{}{$22$}{$11$}{}+\wtile{}{$11$}{}{}+\wtile{}{}{$11$}{}+\wtile{}{}{}{}.
    \end{align*}

    We tile these tiles on a chessboard with only half the width, where the right boundary corresponds to the vertical symmetry axis of the original chessboard.
    Therefore, we no longer restrict the characters on the right boundary to be $\sharp$.
    Instead, we restrict that the two characters at the same position are identical.
    That is, perform the following substitution on the tile set for the rightmost column:
    \[
        \wtile{$\alpha\beta'$}{$aa'$}{$bb'$}{$\beta\alpha'$}\mapsto
        \begin{cases}
            \wtile{$\alpha\beta'$}{$aa'$}{$bb'$}{}, & \beta=\alpha'; \\
            0, & \beta\neq \alpha'.
        \end{cases}
    \]
    For this example, the resulting tile set is
    \[
        \wtile{$12$}{}{}{}+\wtile{$11$}{}{}{}+\wtile{}{}{$22$}{}+\wtile{}{$22$}{$11$}{}+\wtile{}{$11$}{}{}+\wtile{}{}{$11$}{}+\wtile{}{}{}{}.
    \]

    And if the width is odd, we no longer perform the substitution described above.
    Instead, we add a column on the right, consisting of tiles obtained by horizontally folding the tiles from the original tile set, that is:
    \[
        \wtile{$\alpha$}{$a$}{$b$}{$\beta$} \mapsto \wtile{$\alpha\beta$}{$a$}{$b$}{},
    \]
    where \wtile{$\alpha$}{$a$}{$b$}{$\beta$} range over all of those Wang tiles that may appear in positions with left-right symmetry within the polyominoes.
    For this example, the resulting tile set is
    \[
        \wtile{$21$}{}{}{}+\wtile{}{}{$2$}{}+\wtile{}{$2$}{$1$}{}+\wtile{}{$1$}{}{}+\wtile{}{}{$1$}{}+\wtile{}{}{}{}.
    \]

    Such restrictions on symmetry do not significantly increase the computational complexity.
    On the contrary, the computational complexity can likely to be decreased, as the folded tiles may contain fewer pairs of horizontal and vertical characters than the total possible number, while the width is reduced to half or approximately half.
    In this example, the folded tiles contain the following horizontal and vertical character pairs (i.e., the new alphabets for positions excluding the rightmost column)
    \[
        \hat{\varTheta}= \{\sharp,11, 12,21 \}, \hat{\varSigma}= \{\sharp,11,22 \},
    \]
    while their sizes are $\hat{\theta}=4$, $\hat{\sigma}=3$.

    Consequently, for the width $l$, the size of the state tensor is reduced from $\theta\sigma^l$ to
    \[
        \begin{cases}
            \hat{\theta}\hat{\sigma}^{\frac{l}{2}}, & 2\mid l; \\
            \hat{\theta}\sigma\hat{\sigma}^{\frac{l-1}{2}}, & 2\nmid l.
        \end{cases}
    \]
    It is not difficult to see that this reduction is exponential when $\hat{\sigma}<\sigma^2$.

    In addition, one may transpose the chessboard and apply the same method to fold it vertically and perform the enumeration.
    In this case, the size of the tensor (or, more precisely, the size during the majority of the tiling progress) becomes
    \[
        \hat{\sigma}\cdot \hat{\theta}^h.
    \]
    When performing the enumeration, one should choose the approach that results in a smaller tensor size to enhance computational efficiency.

    Based on the enumeration with symmetry restriction, by applying Polya's enumeration theorem, we can also compute the tiling enumeration that flip symmetries are not distinguished.

    \subsection{Open problems}

    As demonstrated in Section 4, the key to transforming the tiling of polyominoes (or their variants) into Wang tiling lies in the Wang decomposition.
    This, however, leaves numerous questions pertaining to Wang decomposition unresolved, such as:

    \begin{problem}
        How can we find the simplest Wang partition (if it exists) for a given set of polyominoes (or their variants)?
    \end{problem}

    \begin{problem}
        Under what conditions can an infinite set of polyominoes (or their variants) admit a Wang partition?
    \end{problem}

    Moreover, the enumeration problem for tilings of what scale remains beyond the computational capacity of our algorithms.
    Consequently, we are compelled to inquire:
    \begin{problem}
        Is there potential to further simplify computations for certain special tiling enumeration problems by leveraging the properties of tensors?
    \end{problem}

    \begin{problem}
        Is it possible to derive specific form of the enumeration formulas (or prove their impossibility) based on the tensor operations we defined?
    \end{problem}

    For tilings of excessively large scales, we can also provide estimations.
    From the computational results, it can be observed that when $l, h\rightarrow \infty$ ,the growth trend of most of the enumeration of tilings within $l \times h$ chessboards approximates to
    \[
        N(l,h) \sim \lambda^{lh},
    \]
    where $\lambda$ is called the \textit{entropy constant}~\cite{entropy1, entropy2}; or, satisfy such approximation periodically.
    For instance, in the case of $\Iii+\Oiv$ tiling enumeration, we can provide the approximation:
    \[
        N(l, h) \sim
        \begin{cases}
            \lambda^{lh}, & 2\mid lh;\\
            0, & 2 \nmid lh,
        \end{cases}
    \]
    and $\lambda\approx1.46$, a rough estimate of the constant based on existing results.

    Hence, we are prompted to inquire:
    \begin{problem}
        Can a universal method be found to provide approximations for large-scale cases of polyomino tiling enumeration?
        How to provide more accurate estimates of the corresponding constants?
    \end{problem}

    \bibliographystyle{unsrt}
    \bibliography{main}

\end{document}